\documentclass [11pt]{amsart}
\usepackage {amsmath, amssymb, amscd, amsthm, mathrsfs, url ,pinlabel,hyperref, verbatim, lipsum}
\usepackage[text={5.5in,9.1in},centering,letterpaper,dvips]{geometry}
\usepackage{color,multirow,dcpic,latexsym,pictexwd,graphics,pdfpages}
\usepackage{fancyhdr, tikz, setspace}
\usetikzlibrary{positioning, arrows.meta}
\usepackage{scrextend}

\usepackage[T1]{fontenc}
\usepackage{lmodern}
\DeclareUrlCommand{\bfurl}{}

\newtheorem {theorem}{Theorem}
\newtheorem {lemma}[theorem]{Lemma}

\newtheorem {question}[theorem]{Question}

\theoremstyle{remark}

\newtheorem {example}[theorem]{Example}

\numberwithin{equation}{section}
\numberwithin{theorem}{section}

\newcommand{\R}{\mathbb{R}}
\newcommand{\C}{\mathbb{C}}
\newcommand{\diag}{\textrm{diag}}
\newcommand{\lk}{\textrm{lk}}
\newcommand{\SymZ}{\mathbb{Z}}

\DeclareMathOperator{\Spin}{{spin}^c}
\DeclareMathOperator{\rank}{rank}
\DeclareMathOperator{\d3}{d_3}
\DeclareMathOperator{\cofac}{cofac}

\title{Real Algebraic Overtwisted Contact Structures on 3-spheres}
\author{\c{S}eyma Karaderel\.{I}}
\address{Department of Mathematics, Bo\u{g}az{\i}\c{c}{\i}  University, Bebek, Istanbul, 34342, Turkey}
\email{\href{mailto:seyma.karadereli@boun.edu.tr}{seyma.karadereli@boun.edu.tr}}
\author{ Fer\.{i}t \"{O}zt\"{u}rk}
\address{Department of Mathematics, Bo\u{g}az{\i}\c{c}{\i}  University, Bebek, Istanbul, 34342, Turkey}
\email{\href{mailto: ferit.ozturk@boun.edu.tr}{ferit.ozturk@boun.edu.tr}}
\date{}

\begin{document}

\begin{abstract}
A real algebraic link in the 3-sphere is defined as the zero locus in the 3-sphere of a real algebraic function from $\R^4$ to $\R^2$.
A real algebraic open book decomposition on the 3-sphere is by definition  the Milnor fibration of such  a real algebraic function, in case it exists. 
We prove that every overtwisted contact structure on the 3-sphere with positive three dimensonal invariant $d_3$ (apart from possibly 9 exceptions)
are real algebraic. Our construction shows in particular that the pages of the associated  open books are planar.
\end{abstract}
\maketitle

\section{Introduction}

A Milnor fillable 3-manifold is a connected closed oriented contact 3-manifold 
which is contact isomorphic to the contact link manifold of a complex analytic surface with isolated singularity. 
We know that any such manifold admits a unique Milnor fillable contact structure up to contactomorphism \cite{CNP} and moreover a Milnor fillable contact structure
is tight. For instance there is a unique tight contact structure on the 3-sphere $S^3$ and it is Milnor fillable (by e.g. the nonsingularity $0$ in $\C^2$).

Here we ask a similar question regarding overtwisted contact structures. 
We confine ourselves to $S^3$ although the definitions and questions below can be easily generalized.
We investigate fibered links in $S^3$ which are given real algebraically (or more generally real analytically).
Let us call an oriented link in $S^3$  weakly real algebraic if it is isotopic to the link of a real algebraic surface with an isolated singularity at 0  (i.e. it is the zero locus of an algebraic map $h:\R^4\rightarrow\R^2$ with an isolated critical point on its zero locus). It is well-known that every  link in $S^3$ is weakly real algebraic \cite{AK}. 
Let us call an oriented link in $S^3$ real algebraic if the map $h$ has an isolated critical point at 0 in $\R^4$. This condition is called the Milnor condition.
In such a case $h$ determines also a Milnor fibration in $S^3$ \cite{M}. In other words the real algebraic link is the binding of an (in general rational) open book 
with the open book decomposition given as the Milnor fibration (see e.g. \cite{BaEt} for rational open books).

Compared to weakly real algebraic setup it is rather hard to construct examples of real algebraic maps with an isolated singularity and this has been long studied. 
For example it is known that figure-eight knot is not complex algebraic but is real algebraic \cite{Perron}. We believe it is still unknown whether there exists
a weakly real algebraic link in $S^3$ which is not real algebraic (see e.g. \cite{Bo}).

An obvious way to produce real algebraic links in $S^3$ is as follows. Take two nonconstant complex algebraic maps 
$f,g:\C^2\rightarrow \C$ and consider the real algebraic map $h=f\overline{g}$. 
The oriented link  $L$  that is the zero locus of $h$ in $S^3$ has components $\{f=0\}\cap S^3$ with canonical orientations and  $\{g=0\}\cap S^3$ with the reverse orientations. It turns out that these links are graph links, i.e. spliced Seifert links \cite{EN}. 
Moreover $h$ has an isolated singularity at 0 if and only if $L$ is fibered, and in that case the Milnor fibration coincides with that determined by $h/|h|$ \cite{Pichon}.
Now, as a corollary to \cite{Ishi} the real algebraic open book corresponding to such $h$ determines an overtwisted contact structure on $S^3$.

Recall that there are countably infinite number of overtwisted structures in $S^3$. They are distinguished by the half-integer-valued $d_3$ invariant (see e.g. \cite{DGS}).
Regarding our question we prove

\begin{theorem}
All overtwisted contact structures on $S^3$ with positive $d_3$ with $d_3+\frac{1}{2}\not\in\{4, 11, 17, 19, 47, 61, 79, 95, 109\}$ are real algebraic.

They are given as the oriented links of real algebraic functions of the form  $f\overline{g}$ for $f,g:\C^2\rightarrow\C$ complex algebraic.
Moreover the associated real algebraic open books have planar pages.
\label{d3pos}
\end{theorem}

We believe that the nine sporadic exceptions that appear in the theorem are real algebraic as well, although  the families of real algebraic Milnor fibrations that we have produced miss them. The nonnegativity that emerge might be more resilient. Thus we ask

\begin{question}
Is there a real algebraic overtwisted contact structure on $S^3$ with negative $d_3$? If yes, the supporting real algebraic open book is rational in general.
Can it be an honest open book? Can it have planar pages?
\end{question}

Generalizing our definitions we ask 
\begin{question}
Is it true that every overtwisted contact structure on a Milnor fillable 3-manifold is real algebraic?
\end{question}

To proceed towards the proof of Theorem~\ref{d3pos}, we recall in Section~\ref{sec:multsplice} the Seifert and graph multilinks and the splicing operation. There we also give
our choice of examples for families of  fibered graph multilinks in $S^3$ and compute the associated monodromy maps. 
In Section~\ref{sec:alg} we demonstrate that those families of graph multilinks and the corresponding open book decompositions are real algebraic.
In Section~\ref{sec:4mfd}  we briefly recall a way to compute the $d_3$ invariant, by constructing almost complex 4-manifolds that fill the given open book decompositions in $S^3$. 
Finally in Section~\ref{mainproof} we prove Theorem~\ref{d3pos} by computing the $d_3$ invariants explicitly for our 
proposed families of examples. 
It turns out that one of our families of graph multilinks exhausts all the overtwisted structures with $d_3>461$. Then by  computer aid we 
show that those with $0<d_3<461$ (except the 9 values given in the theorem) are realized by our families of graph multilinks as well.
In the computation of $d_3$ the constructed 4-manifolds have large  intersection matrices. 
For the clarity of the exposition, we present those intersection matrices 
in Appendix~\ref{thematrices} and the tedious computations regarding these matrices are given in Appendix~\ref{theinvariants}.

\section{Seifert multilinks and Splicing }
\label{sec:multsplice}
In this section we recall introductory information on Seifert and graph multilinks and present several families of examples.
Our discussion is based on \cite{EN}.
\subsection{Seifert multilinks}  
\label{sec:multilink}
A Seifert fibered manifold is a closed 3-manifold given as an $S^1$-bundle with the orbit space a 2-orbifold. 
A Seifert multilink in a Seifert fibered 3-sphere is an oriented link $L$ that is constituted of a finite number of Seifert fibers $S_i$ 
and an integer  multiplicity $m_i$ assigned to each component. In this work we are solely interested in Seifert multilinks in $S^3$.
We are going to denote a Seifert multilink with $n$ components by $L(m_1, . . . , m_n)$. 
Note that $L$ is canonically oriented by the sign of the multiplicities $m_i$. 
In this setup the  homology class $\underline{m}=(m_1,\cdots,m_n)\in H_1(L)\simeq \mathbb{Z}^n$  
determines a cohomology class in the link complement as well, since $H_1(L)\simeq H^1(M-L)$.  That class is given by 
\begin{align*}
 \underline{m}(\gamma)=\lk(L, \gamma)=\sum_{i=1}^n m_i \cdot \lk(S_i, \gamma).
\end{align*}
Let $\mu_i$ denote the meridian of the $i$-th link component. Then we have $\underline{m}(\mu_i)=m_i$.  Moreover we can realize the Seifert surface of the  multilink 
as an embedded oriented surface whose intersection with the boundary of a tubular neighborhood of $S_i$ is $\left(\delta_i\cdot (\frac{m_i}{\delta_i},\frac{-(m_i)^\prime}{\delta_i})\right)$-cable of $S_i$, where $(m_i)^\prime=\sum_{j\neq i}^n m_j \lk(S_i,S_j)$ and $\delta_i=\gcd(m_i, m_i^{\prime})$ \cite[p.30]{EN}. 

\begin{figure}[h]
	\begin{center}
\resizebox{13cm}{!}{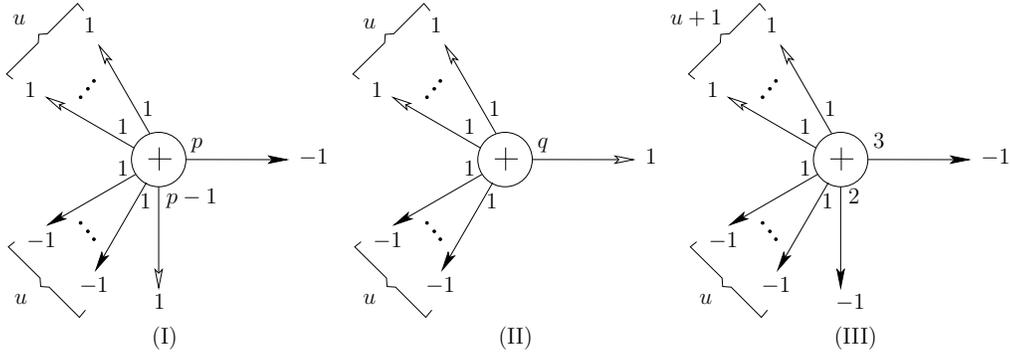}
\caption{Splice diagrams for Seifert multilinks of type (I), (II), (III)}
		\label{fig:SpliceDiagrams}
	\end{center}
\end{figure}

Multilinks are represented by splice diagrams as illustrated in Figure~\ref{fig:SpliceDiagrams}.  
The central node represents the ambient Seifert manifold. The numbers adjacent to the node for each branch are called the weights and the numbers next to the 
arrowheads are the multiplicities above.

An arrowhead vertex with weight $+1$ (respectively $> 1$) corresponds to a  regular (respectively singular) Seifert fiber.  
The multilink (I) in Figure~\ref{fig:SpliceDiagrams} has $2u+2$ connected components in the underlying manifold $S^3$ on which the Seifert fibration is given 
by the $S^1$-action $(x,y)\mapsto (t^{p-1}x,t^{p}y)$  for $t\in S^1$.  
Here  the orbit $\{x=0\}$ corresponds to the singular fiber $S_1$ with weight $p$ and $\{y=0\}$ corresponds to the singular fiber $S_2$ with weight $p-1$. 
The linking numbers of link components can be computed easily using the splice diagram  \cite[Proposition~7.4]{EN}.   The linking number of any nonsingular fiber with the singular fiber $S_1$ and respectively with $S_2$  is the product of weights of the remaining vertices, which equals $p-1$ and $p$ respectively.  
The linking number of $S_1$ and $S_2$ is $1$. Thus  the multilink (I) is isotopic to the negative Hopf link with $u$ positively oriented and $u$ negatively oriented copies of the $(p,p-1)$ torus knot cabled around $S_1$.  

A multilink $L(\underline{m})$ is  fibered if there exists a locally trivial fibration  $M-L\to S^1$ in the homotopy class corresponding to $\underline{m}$,  whose fibers are minimal Seifert surfaces for the multilink.  
Using the analytic description of the Seifert fibration of the link exterior,  it can be easily seen that a Seifert multilink is fibered if and only if the linking number of any nonsingular fiber $\gamma$ with the multilink does not vanish \cite[p.90]{EN}.  In other words,  denoting  the weight of the $i$-th link component $S_i$ as $\alpha_i$,
$$l=\underline{m}(\gamma)=\sum_{i=1}^n m_i \lk(\gamma,S_i)=\sum_{i=1}^n m_i\alpha_1\cdots\hat{\alpha_i}\cdots\alpha_n\neq 0.$$
A fibered multilink determines a rational open book decomposition for the ambient Seifert manifold.
If each $m_i=\pm 1$ then the open book is an integral open book. 
We recall that if $l=1$ the pages of the open book are planar.

The monodromy of the fibration can be represented as the flow along the Seifert fibers.  
Thus  in the interior of the pages it is isotopic to a homeomorphism of order $l$. On the other hand the monodromy flow near the boundary components is 
computed as a $\displaystyle (-\frac{\delta_i}{m_i l}\alpha_i)$-worth (in general rational) twist along a boundary parallel curve \cite[p.108]{EN}.

\begin{example}
For the multilinks of type (I) given in Figure~\ref{fig:SpliceDiagrams},  the multilink is fibered since we have  $l=-1\cdot (p-1)+1\cdot p+u\cdot (1)\cdot p(p-1)+u\cdot(-1)p(p-1)=1\neq 0$.  The pages are  $(2u+2)$-punctured spheres. 
We recall that the families of diagrams in Figure~\ref{fig:SpliceDiagrams} are exactly the three of four types of Seifert multilinks with $l=1$ \cite[p.123]{EN}
The monodromy flow is trivial in the interior of the pages.  However near the boundary components corresponding to the singular fibers, the flow is given as 
$-\frac{1}{-1\cdot 1}p=p$ and $-\frac{1}{1\cdot 1}(p-1)=-(p-1)$ twists.    Along the  boundary components corresponding to the nonsingular fibers with positive and negative multiplicities, the flow is  $-1$ and  $+1$ twist respectively.  Therefore the monodromy is given as 
\begin{equation}
\label{monodromy(I)}
\phi=a^{p}\cdot b^{-(p-1)}\cdot c_{1}^{-1} \cdots c_{u}^{-1} \cdot d_{1} ^{1}\cdots d_{u}^{1}.
\end{equation} 
Here, $a$  and $b$ denote Dehn twists along curves parallel to the boundary components $\{x=0\}$ and  $\{y=0\}$ respectively; $c_{i}$ and $d_i$ are twists along curves parallel to the nonsingular components with positive and negative multiplicities respectively.

Similarly, as noted above, the multilinks of type (II) and type (III) in Figure~\ref{fig:SpliceDiagrams} are fibered multilinks  in $S^3$ with $l=1$ too.  
The pages of the multilink of type (II) are $(2u+1)$-punctured spheres and the monodromy is
\begin{equation}
\label{monodromy(II)}
\phi=a^{-q}\cdot b ^{-1}\cdot c_{1}^{-1}\cdots c_{u-1}^{-1} \cdot d_{1}^{1} \cdots d_{u}^{1}.
\end{equation}  
The pages of the multilink of type (III) are $(2u+3)$-punctured spheres and the monodromy is
\begin{equation}
\label{monodromy(III)}
\phi=a^{3}\cdot b^{2}\cdot c_{1}^{-1} \cdots c_{u+1}^{-1} \cdot d_{1} ^{1}\cdots d_{u}^{1}.
\end{equation}
\end{example}

\subsection{Splicing multilinks}
\label{sec:splicing}
 The splice of two multilinks along a specified pair of  link components is constructed topologically
by excising tubular neighborhoods of the given link components and gluing the remaining manifolds in a meridian-to-longitude fashion.  
Note that topologically splicing multilinks in $S^3$ produces still a multilink in $S^3$.
Moreover a cohomology class is determined by the multiplicities of the components of the resulting multilink.
For the splicing operation we also require that the restriction of this cohomology class on each manifold gives the cohomology class of the splice component. 
This condition is equivalent to the following. Let $S_0$ and $\tilde{S_0}$ with multiplicities $m_0$ and $\tilde{m_0}$ be the spliced link components;  
$({\mu_0},{\lambda_0}),(\tilde{\mu_0},\tilde{\lambda_0})$ be the meridians and longitudes on the tori on which the splicing occurs. Then we have
\begin{equation}
\begin{split}
&m_0=  \underline{m}(\mu_0)= \underline{\tilde{m}}(\tilde{\lambda}_0)=(\tilde{m}_0)^\prime, \\
&\tilde{m}_0= \underline{\tilde{m}}(\tilde{\mu}_0)= \underline{m}(\lambda_0)=(m_0)^\prime.
\end{split}
\label{eq:SpliceProperty}
\end{equation}
where $(\tilde{m}_0)^\prime, (m_0)^\prime$ are defined as in Section~\ref{sec:multilink}. Observe that  these requirements are exactly the conditions for the Seifert surfaces in each splice component to glue together along the splicing tori.
Moreover  since Seifert surfaces approach the spliced link components as $\delta_0=\gcd(m_0, (m_0)^\prime)$ copies of the 
$(\frac{m_0}{\delta_0}, \frac{(m_0)^\prime}{\delta_0})$ curve,  the Seifert surfaces are pasted together along $\delta_0$  tori.

As an example,  consider the multilink (I) in Figure~\ref{fig:SpliceDiagrams} and there  the link component $S_1$ of weight $p$.  
Since $\underline{m}(\lambda_1)=(1)\cdot 1\cdots 1+u( (1)\cdot 1 \cdots 1\cdot(p-1))+u( (-1)\cdot1\cdots 1\cdot (p-1))=1$ 
and $\underline{m}(\mu_1)=-1$, one can splice $S_1$ only with a link component whose multiplicity is $\tilde{m}_1=1$ and 
$(\tilde{m}_1)^\prime=-1$, i.e.  the pages must approach the link component as $(1,-1)$ curves.  Similarly for the link component $S_2$ of weight $p-1$, we have 
$\underline{m}(\lambda_2)=(-1)\cdot 1\cdots 1+u((1)\cdot1\cdots 1\cdot p)+ u((-1)\cdot1\cdots 1\cdot p)=-1$ and 
$\underline{m}(\mu_2)=1$. Therefore,  given two multilinks of type (I)  one can only splice the former link component 
of one  with the latter link component of the other.

Another possible splicing occurs between the splice multilink (II) with $q=2$ and the  multilink (I). 
Since the multiplicity of the link component with weight $q$ of the multilink (II) is 1,  splicing it with  the link component with weight $p$ of the multilink (I) is possible if  $\underline{m}(\lambda)=(u-1)\cdot 1\cdot q +u\cdot (-1)\cdot q+1 = 1-q = -1$ for the singular fiber of the multilink (II).  

\begin{figure}[h]
	\begin{center}
\resizebox{13cm}{!}{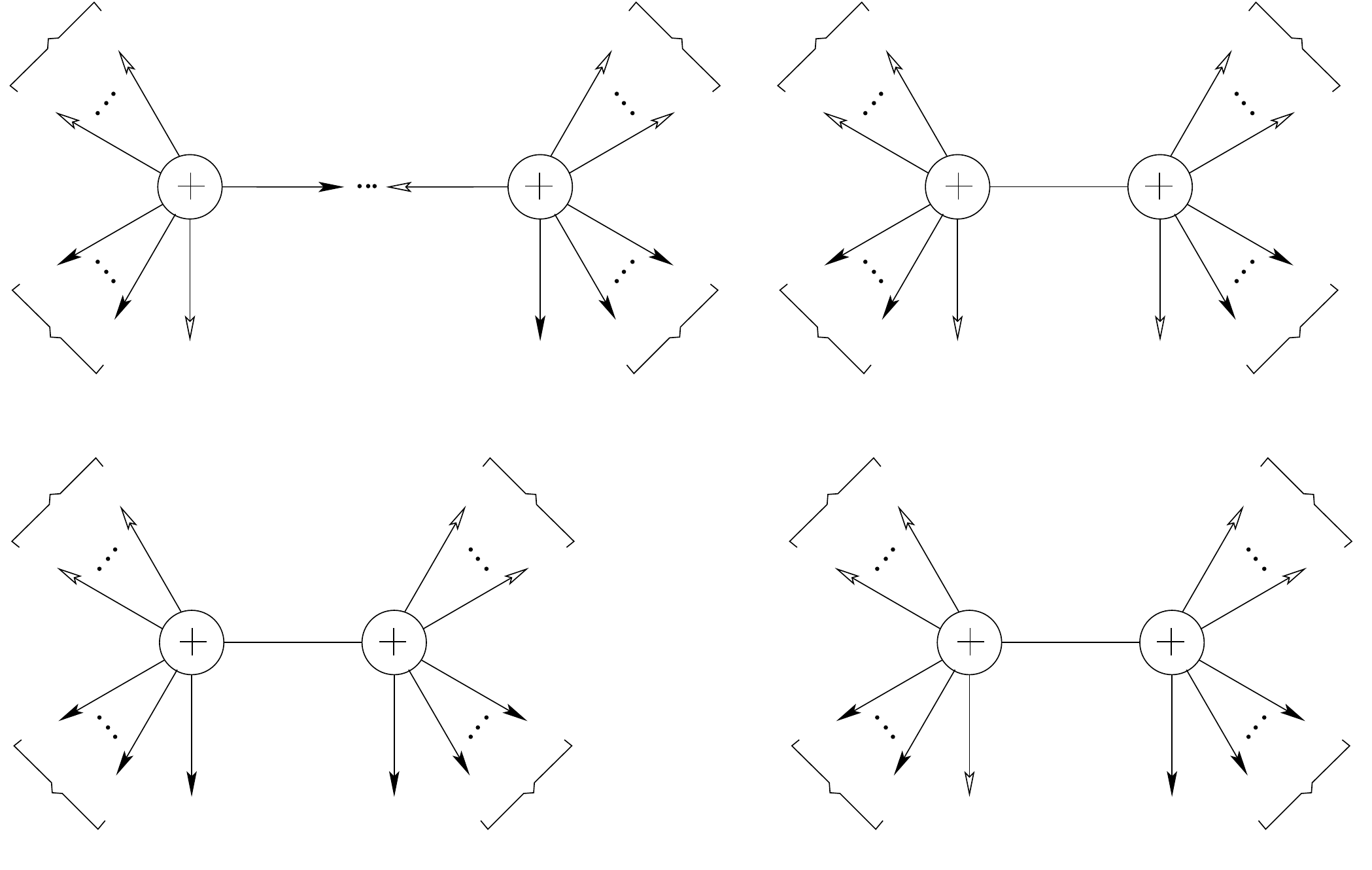}
		\caption{Splice diagrams for (a): I-I, (b): II-I, (c): III-I, (d): II-III.}
		\label{fig:AllSplice}
	\end{center}
\end{figure}

Going through all possible cases we obtain the following list. 
\begin{lemma}
All the possible splicing operations between the multilinks (I),(II), (III)  in $S^3$ are
\begin{itemize}
\item $(\textrm{I-I-I-}\ldots)$: here each splicing is as in Figure~\ref{fig:AllSplice}(a). 
\item $(\textrm{II-I-I-}\ldots)$: the first splicing is as in Figure~\ref{fig:AllSplice}(b).
\item $(\textrm{III-I-I-}\ldots)$: the first splicing is as in Figure~\ref{fig:AllSplice}(c).
\item (II-III): the splicing is as in Figure~\ref{fig:AllSplice}(d).
\end{itemize}
\end{lemma}
Splicing of two multilinks is represented by a splice diagram (with more than one node) obtained by joining the two diagrams along the arrowheads corresponding 
to the link components at which splicing occurs. A multilink with such a splice diagram is called a graph multilink.


A graph multilink is fibered if and only if it is an irreducible link and each of its splice components is fibered \cite[Theorem~4.2]{EN}. 
The monodromy is pieced together from the monodromy maps of the splice components.  
In each splice component the monodromy is given by the flow along the corresponding Seifert fibers whereas on the tubular neighborhoods of the separating tori, it has two different flows in each end given by  the Seifert fibration of each Seifert component.   Therefore  after splicing, the Dehn twists corresponding to glued boundaries become 
trivial and on the separating tori  the monodromy acts as a twist map  which measures the difference between the two flows of Seifert fibers.  
In \cite[Theorem~13.1]{EN} the monodromy flow on a separating torus is computed as a $\tau$-worth twist with 
 \begin{equation}
\label{eq:twistsep}
\tau=\frac{-\delta_0}{l_1\cdot l_2}(\alpha_0\beta_0-\alpha_1\cdot\ldots\cdot\alpha_n\cdot\beta_1\cdot\ldots\cdot\beta_m)
\end{equation}
where $\alpha_0$, $\beta_0$ are the weights of the spliced components and  $\alpha_i, \beta_j$ are the weights of the remaining link components.   

\begin{example}
\label{ex:splicingI-I}
Consider the multilink (I-I) given in Figure~\ref{fig:AllSplice}.  Note that when $q=p$  the graph multilink is simply a Seifert multilink \cite[Theorem~8.1(6)]{EN}.  
Then let us consider the case $q>p$.

By the previous discussion we know that $l_1=l_2=1$; also  $\delta=\gcd(-1,1)=1$. Thus $$\tau=-\dfrac{1}{1\cdot 1}(p(q-1)-q(p-1))=p-q.$$

Since $\delta=1$,  we glue the pages of the spliced components, which are  $(2u+2)$- and $(2v+2)$-punctured spheres respectively,  along a single annulus neighborhood of the spliced boundary components. Consequently the pages of the spliced multilink are $(2u+2v+2)$-punctured spheres.

As given in (\ref{monodromy(I)}) the  splice components have monodromies $\phi_1=\alpha^{p}\cdot a^{-(p-1)}\cdot c_{1}^{-1} \dots c_{u}^{-1} \cdot d_{1} ^{1}\dots d_{u}^{1}$ and $\phi_2=b^{q}\cdot \beta^{-(q-1)}\cdot e_{1}^{-1} \dots e_{v}^{-1} \cdot f_{1} ^{1}\dots f_{v}^{1}$.  The  monodromy flow is  $q-p$ negative Dehn
twists about the core circle, say  $\gamma$, in the annulus. Therefore the monodromy of the spliced multilink is
\begin{equation}
\label{monodromy(I-I)}
\phi=a^{-(p-1)}\cdot c_1^{-1}\cdots c_u^{-1} \cdot d_1^{1}\cdots d_u^{1} \cdot  \gamma^{-(q-p)}\cdot b^{q} \cdot  e_1^{-1}\cdots e_v^{-1} \cdot f_1^{1} \cdots f_v^{1};
\end{equation}
\end{example}

\begin{example}
Similarly let us consider a graph multilink of the form (I-I-I) as  in Figure~\ref{fig:Splice111}.
Note that we splice the link component with weight $q$ of the first splice component to the link component with weight $(r-1)$ of the second splice component.

\begin{figure}[htbp]
	\begin{center}
\resizebox{9cm}{!}{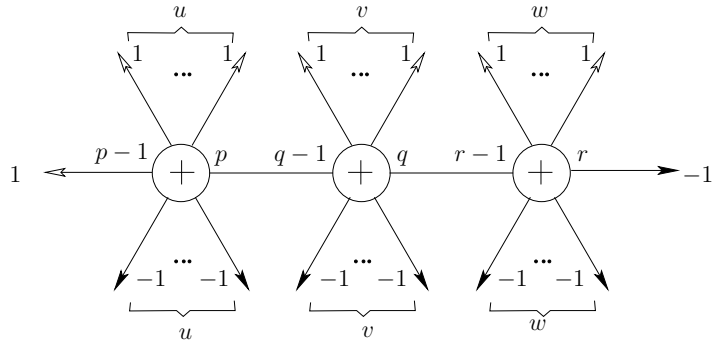}
		\caption{Splice diagram for (I-I-I).}
		\label{fig:Splice111}
	\end{center}
\end{figure}
As in the previous examples  $l_1=l_2=1$ and $\delta=\gcd(m_1,m_2)=1$. Assuming $r >q$,  we have $$\tau=-\dfrac{\delta}{l_1 l_2}(q(r-1)-r(q-1)=q-r < 0.$$
The page of the splice multilink is a union of the pages of the splice components joined together along a boundary by a  $(q-r)$-twisted  annulus (since $\delta=1$). 
Since the splice components have $(2u+2v+2)$- and $(2w+2)$-punctured sphere pages, the pages for the splice link are  $(2u+2v+2w+2)$-punctured spheres.

The monodromy of the new fibration is 
\begin{equation}
\label{monodromy(I-I-I)}
\phi= a^{-(p-1)} c_1^{-1}\cdots c_u^{-1}  d_1^{1}\cdots d_u^{1}   \gamma^{-(q-p)} e_1^{-1}\cdots e_v^{-1}  f_1^{1} \cdots f_v^{1} \theta^{-(r-q)}  b^{r} g_1^{-1}\cdots g_w^{-1}  h_1^{1} \cdots h_w^{1}
\end{equation} 
where $\theta$ denotes the Dehn twist about the core  circle in the  latter annulus.
\end{example}

\begin{example}
As in the previous example one can compute the monodromies of the other multilinks given in Figure~\ref{fig:AllSplice}; the monodromy of the splicing (II-I) is 
\begin{equation}
\label{monodromy(II-I)}
\phi=a^{-(p-1)}\cdot c_1^{-1}\cdots c_u^{-1}\cdot d_1^{1}\cdots d_u^{1}\cdot \gamma^{(p-2)}\cdot b^{-2} \cdot e_1^{-1}\cdots e_{v-1}^{-1}  \cdot f_1^{1} \cdots f_{v}^{1}
\end{equation}
Here,  we assume that $p\geq 3$  because  the graph multilink is simply a Seifert multilink when $p=2$ \cite[Theorem~8.1(6)]{EN}.

The monodromy of the splicing (III-I) is 
\begin{equation}
\label{monodromy(III-I)}
\phi=a^{(p)}\cdot c_1^{-1}\cdots c_u^{-1}\cdot d_1^{1}\cdots d_u^{1}\cdot \gamma^{-(p-3)}\cdot b^{2} \cdot e_1^{-1}\cdots e_{v+1}^{-1}  \cdot f_1^{1} \cdots f_{v}^{1}
\end{equation} 
For the same reason as before,  we assume $p\geq 4$ here.

The monodromy of the splicing (II-III) is 
\begin{equation}
\label{monodromy(II-III)}
\phi=a^{-2}\cdot c_1^{-1}\cdots c_{u-1}^{-1}\cdot d_1^{1}\cdots d_u^{1}\cdot \gamma^{1}\cdot b^{2} \cdot e_1^{-1}\cdots e_{v+1}^{-1}  \cdot f_1^{1} \cdots f_{v}^{1}
\end{equation}
\end{example}

\section{Real algebraic singularities and associated contact structures}
\label{sec:alg}

In this section we assert that the graph multilinks and the associated open books that have been considered in the previous section can be realized via real algebraic singularities.


For an isolated singularity of a holomorphic (or a complex algebraic) function from $\mathbb{C}^2$ to $\mathbb{C}$, 
the corresponding Milnor fibration defines an open book structure on $S^3$, whose binding is isotopic to the singularity link. In such a setup we call the singularity link, 
the open book and the supported tight contact structure complex analytic/algebraic.
Any complex algebraic link in $S^3$  is a graph multilink and the corresponding splice diagram can be deduced from the Puiseux pairs \cite[Appendix~1]{EN}. 
Of course not all the graph multilinks in $S^3$ are complex algebraic. \cite[Theorem~9.4]{EN} gives the precise condition for a graph multilink to be complex algebraic.

Similarly an isolated singularity of a  real analytic function $f :\mathbb{R}^4 \to \mathbb{R}^2$ determines a Milnor fibration in $S^3$ 
under the condition that the Jacobian matrix of $f$ has rank $2$ on an open neighborhood of the origin, except the origin. 
This  is  the Milnor condition. A link is  said to be real analytic/algebraic  if it is the singularity link of a real analytic/algebraic map 
$f :\mathbb{R}^4 \to \mathbb{R}^2$ that satisfies the Milnor condition.
In the absence of the Milnor condition,  $f$ might not even have an isolated singularity and a Milnor fibration.    
In \cite{Pichon} and  \cite{PS},  the authors examine real algebraic germs with isolated singularities of mixed functions  of the form $f\overline{g}$ where $f$ and $g$ are holomorphic functions.  They discuss the Milnor fibration in the link exterior and   the  geometry of the fibration near the singularity link.

The isotopy class of a multilink is encoded in a plumbing tree that is decorated with arrows having multiplicities for the link components.  
When a multilink is isotopic to the singularity of a holomorphic germ,  the plumbing tree for the multilink can be obtained 
as the dual tree of any normal crossing resolution of the function.  Since $L_{f\overline{g}}$ as an unoriented link is $L_f\cup L_g$, it can be seen that the resolution graph of a real algebraic germ of the form $f\bar{g}$ is nothing but the resolution graph of $fg$ with negative signs for the multiplicities of the link components corresponding to $g$.  The following theorem explains when the singularity link of a real algebraic germ of the form $f\bar{g}$ has a Milnor fibration.
\begin{theorem}\emph{\cite[Theorem~5.1]{Pichon}}
\label{thm:realalgebraic}
Let $f:\left(\mathbf{C}^{2}, 0\right) \rightarrow(\mathbf{C}, 0)$ and $g:\left(\mathbf{C}^{2}, 0\right) \rightarrow(\mathbf{C}, 0)$ be two holomorphic germs with isolated singularities and having no common branches. Then the real analytic germ $f \bar{g}$ has an isolated singularity at 0 if and only if the link $L_{f}-L_{g}$ is fibered.

Moreover, if this condition holds, then the associated Milnor fibration is a fibration of the link $L_{f}-L_{g}$.
\end{theorem}


\begin{figure}[htbp]
	\begin{center}
		\includegraphics[width=0.8\columnwidth]{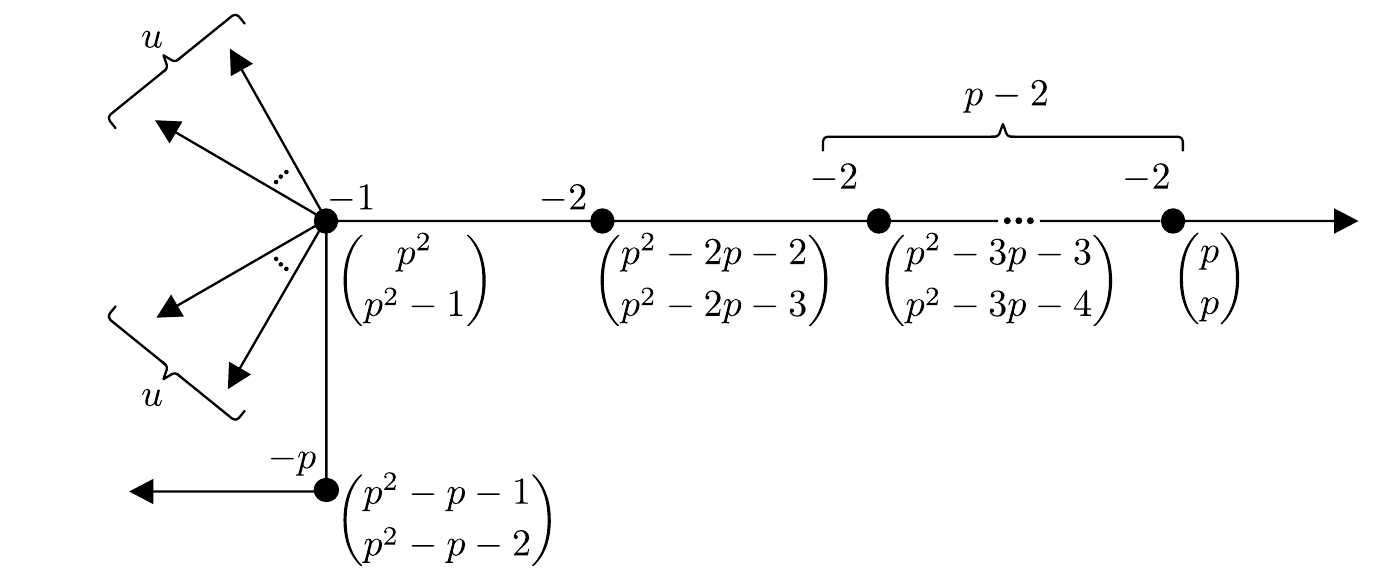}
		\caption{Dual tree of a resolution $\pi$ of $fg$ with associated   multiplicities given in the parantheses which are the multiplicities  $m_i^f$, $m_i^g$ of $f\circ \pi$ and $g\circ \pi$  respectively along the irreducible component for the $i$-th exceptional divisor.  As a side remark we recall that $L_{f}-L_{g}$ is fibered if and only if $m_i^f\neq m_i^g$ at the rupture vertices \cite[Corollary~2.2]{PS}. }
		\label{fig:Plumbing}
	\end{center}
\end{figure}
\begin{example} 
\label{ex:real}
For  an $\eta^{2u+1}=1$ consider the functions 
\begin{equation}
\label{realalg(I)}
f(x,y)=y  \prod_{i=1}^{u}(x^{p} + \eta^{i} y^{p-1})\text{ and }g(x,y)=x \prod_{j=u+1}^{2u}(x^{p} + \eta^j y^{p-1})
\end{equation}  

After resolving the germ of $f g$,  we obtain the plumbing diagram of $L_{f\bar g}$  given in Figure~\ref{fig:Plumbing}. 
As in \cite[Section~20]{EN},  we can obtain the splice diagram of the singularity link from the plumbing diagram and see that it is isotopic to the multilink of type (I) 
in Figure~\ref{fig:SpliceDiagrams}.  Since we have already noted that the multilink is fibered,  it follows from Theorem~\ref{thm:realalgebraic} 
that $f\bar{g}$ has an isolated singularity and the fibration of the multilink which we investigate in the previous section is the Milnor fibration of the germ.  
Observe also  that the branch  $\{\bar{x}=0\}$ corresponds  to the singular link component of weight $p$,  $\{y=0\}$ corresponds to the singular component of weight $p-1$ 
and the positively (respectively negatively) oriented $u$ copies of $(p,p-1)$ cables around $\{x=0\}$ component correspond 
to the branches $\{ \prod_{i=1}^{u}(x^{p} + \eta^{i} y^{p-1})=0\}$  (respectively $\{ \prod_{i=1}^{u}\overline{(x^{p} + \eta^{i} y^{p-1})}=0\}$).
\end{example}

\begin{example}
Similarly we observe that the singularity links of the  real algebraic germs 
\begin{equation}
\label{realalg(II)}
\displaystyle \left( xy\prod_{i=1}^u(x^q+\eta^{i}y)\right) \cdot \left(\prod_{j=1}^{u-1}\overline{(x^q+\eta^{u+j}y)} \right)
\end{equation} and 
\begin{equation}
\label{realalg(III)}
\left(\prod_{i=1}^{u+1}(x^3+\eta^{i}y^2)\right) \cdot  
\left( \overline{x} \overline{y} \prod_{j=1}^{u}\overline{(x^3+\eta^{u+j+1}y^2)}\right)
\end{equation} 
are isotopic to the fibered multilinks of type (II) and (III) in Figure~\ref{fig:SpliceDiagrams} respectively,  therefore have isolated singularities at the origin and engender  Milnor fibrations. 
\end{example}

As for the graph multilinks obtained via splicing in the previous section, 
a priori they might not be algebraic.  
Consider the canonically oriented graph multilink  isotopic to the multilink (I-I).  This multilink is complex algebraic when $q>p$,  
thanks to the condition given in \cite[Theorem~9.4]{EN}.  The corresponding holomorphic function can be easily deduced 
from the holomorphic germs related to the spliced components as follows.
Recall that we splice the component corresponding to the branch $\{x=0\}$ of  a multilink $L_1$ of type  (I) with weights for singular fibers $p, p-1$ 
with the component $\{y=0\}$ of a multilink $L_2$ of type (I) with weights $q, q-1$.   
By isotopy, the nonsingular link components 
of $L_1$  which are  $(p,p-1)$ cables of $\{x=0\}$ can be realized as $(p-1,p)$ cables of the $\{y=0\}$ component of $L_1$.  
As we splice,  we remove the spliced link components and keep the remaining ones. The resulting multilink is a positive Hopf link with 
$2u$ many $(p-1,p)$ cables around the  link component $\{y=0\}$ (coming from $L_1$) and $2v$ many 
 $(q,q-1)$ cables around the  link component $\{x=0\}$ (coming from $L_2$).  
Again by isotopy,   $(p-1,p)$ cables around the former component can be seen as $(p,p-1)$ cable around the latter.  
The resulting multilink is  the union of all components of the spliced multilinks except the ones we spliced.  
Thus  the corresponding holomorphic function is nothing but the  product of the algebraic functions corresponding to branches.  
Since the spliced multilink (I-I) is the above multilink where some of the link components are oriented negatively,  
it is still algebraic when $q>p$ and the corresponding real algebraic map is the map where we take the conjugate of the 
algebraic functions corresponding to the branches that are oriented negatively.  The real algebraic map corresponding to this graph multilink is  of the form $f\bar{g}$ and given by 
\begin{equation}
\bar{x} y  \prod_{i=1}^{u}(x^p + \eta^{i} y^{p-1}) \prod_{j=u+1}^{2u}\overline{(x^p + \eta^j y^{p-1})}
\prod_{i=1}^{v}(x^q + \eta^{i} y^{q-1}) \prod_{j=u+1}^{2v}\overline{(x^q + \eta^j y^{q-1})}.
\end{equation}

Similarly,  the graph multilink (I-I-I) is real algebraic when $p<q<r$. The multilinks (II-I) and (III-I) are  real algebraic when $p>2$ and $p>3$ respectively. 
Thus  all the families of multilinks we have constructed are real algebraic.

 In \cite{Ishi} it is proven that if the link components of a fibered multilink in a homology 3-sphere are canonically oriented (or all those orientations are reversed),  then  the multilink is the binding of an open book which supports a tight contact structure; otherwise the supported contact structure is  overtwisted. So  one can conclude 
that the Milnor open books of the real algebraic links we have constructed so far support overtwisted contact structures in $S^3$.
 
\section{Calculation of the 3-dimensional invariant from open books}
\label{sec:4mfd}

In this section  we recall how to detect  the overtwisted contact structures compatible with 
the  Milnor fibered multilinks  constructed in the previous sections using  the monodromy data.

Recall that two overtwisted contact structures on $S^3$ are contact isotopic if and only if they are homotopic as 2-plane fields \cite{Elia2}.  Moreover 
the homotopy class of a 2-plane field is determined by the induced $\Spin$ structure and the $\d3$ invariant (see \cite{Tura},  \cite{G}).
Since $S^3$ has a unique $\Spin$ structure,  the overtwisted structures on $S^3$ are classified by their $\d3$ invariants, 
which take values in $\mathbb{Z}+\frac{1}{2}$  (see e.g.      \cite{DGS}). 
There may be various ways to compute the $\d3$ invariant of a given contact structure. Here we will use the method 
in \cite{EO} to calculate that from the monodromy data of the compatible open book.  

It is known that given an achiral Lefschetz fibration on  a 4-manifold $W$ with fibers $F$ with boundary,  $W$ can be described as   $F \times D^2$ with 2-handles attached 
to some vanishing cycles $\gamma_i$ with appropriate framings.  The Lefschetz fibration on  $W$  induces an open book decomposition and hence a contact structure on  $\partial W$.  
The contact structure induced on $\partial W$ is 
obtained by contact $(+ 1)/( -1)$-surgeries on the Legendrian realizations of the vanishing cycles of respectively negative/positive critical points, each embedded in distinct fibers of the open book; the contribution to the monodromy is respectively a left/right handed Dehn twist about the vanishing cycle. 
In the reverse direction given  a 3-manifold with an open book decomposition, the monodromy data determines 
an achiral Lefschetz fibration on a 4-manifold which on the boundary gives the given open book.

It should be noted that  2-handle attachments with $(-1)$ framing  result in an honest Lefschetz fibration 
carrying  a natural almost complex structure which is the extension of the one on $D^2\times F$.  However,   attaching a 
2-handle with  $(+1)$ framing gives an achiral Lefschetz fibration which does not have a natural almost complex structure that comes from extending the older one.  It is shown in \cite{DGS} that  if $W_0$ is the handlebody decomposition of the 4-manifold admitting the Lefschetz fibration constructed via $k$  $(+1)$-surgeries, 
$W=W_0\# k\mathbb{C}P^2$ (with the same boundary) has a natural  almost complex structure.  When the second cohomology has no torsion 
(where $W$ is assumed to have no 1-handles)
one has the following formula (see \cite{EF} or \cite{EO}) which is the generalization of the similar statement in \cite{DGS}:
 \begin{equation}
 \label{eqn:d3invariant}
\d3(\xi)=\frac{1}{4}(c^2(W)-2\chi(W)-3\sigma(W))+k.
\end{equation}
Here  $\sigma(W)$ and $\chi(W)$ are the signature and the Euler characteristic of $W$. The Chern class $c\in H^2(W;\mathbb{Z})$  
is the Poincar\'e dual to $\sum_{i=1}^n r(\gamma_i)C_i$ where $C_i$ is the cocore of the 2-handle attached along the vanishing cycle $\gamma_i$,
and $r(\gamma_i)$ is the rotation number of $\gamma_i$.  
Since $c(W)|_{\partial W}=c(\xi)$ is zero,  $c(W)\in H^2(W)$ comes from a class in $H^2(W,\partial W)$ thus can be squared.   
A way to calculate $r(\gamma_i)$ on a page is explained in \cite{EO} in detail.  
The rotation number  is  equal to the winding number of the projection of the curve to a page with respect to the orientation on the Kirby diagram obtained by the usual orientation of $D^2$ extended over 1-handles.

\section{Proof of Theorem~\ref{d3pos}}
\label{mainproof}
We have seen that the multilinks in Figure~\ref{fig:SpliceDiagrams} are fibered with planar pages  (see Section~\ref{sec:multilink}) and are real algebraic  (see Section~\ref{sec:alg}).  These multilinks can be spliced  together to build wider families of fibered, real algebraic multilinks  (see Section~\ref{sec:splicing}).  Moreover the corresponding contact structures are overtwisted  
 (see Section~\ref{sec:alg}). 
 
 In this section,  we  calculate their $\d3$ invariants and show what 
overtwisted contact structures on $S^3$ are supported by those real algebraic planar open books. 
This discussion will be tied in Section~\ref{ohfinally} to prove Theorem~\ref{d3pos}.

\subsection{Overtwisted structures via (I)}
\label{ex:(I)} 
We first consider the family of multilinks of type (I).
Recall that  the open books that they determine have pages 
$(2u+2)$ times punctured spheres (denoted by $\Sigma_{0,2u+2}$).
Moreover the monodromy (\ref{monodromy(I)}) of the open book is
$$\phi=a^{p}\cdot b^{-(p-1)}\cdot c_{1}^{-1} \dots c_{u}^{-1} \cdot d_{1} ^{1}\dots d_{u}^{1},$$
where $a,b,c$ are boundary parallel curves. Observe that the number of negative Dehn twists in this expression is $p+u-1$.

\begin{figure}[htbp]
	\begin{center}
		\includegraphics[width=1\columnwidth]{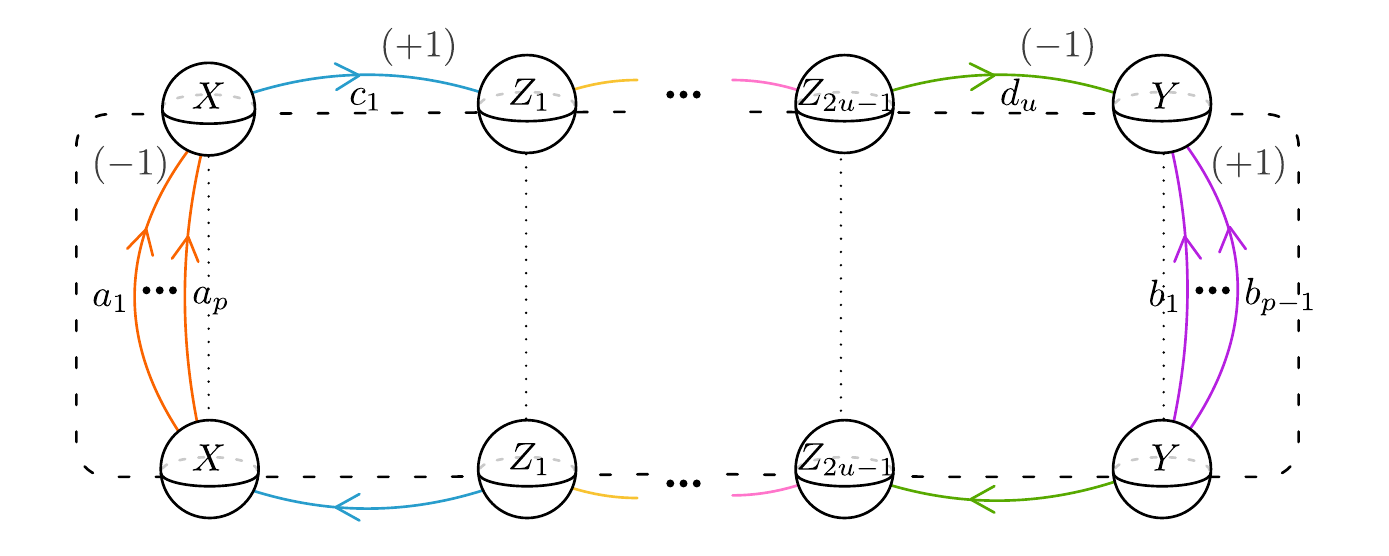}
		\vskip\baselineskip 
		\caption{Kirby diagram for the 4-manifold corresponding to (I)}
		\label{fig:Kirby1}
	\end{center}
\end{figure}

As we discussed in Section~\ref{sec:4mfd}, via the monodromy information of the given open book decomposition we can construct a 4-manifold with boundary $S^3$ as the underlying space of an achiral Lefschetz fibration. In that way we can calculate the $\d3$ invariant of the overtwisted contact structure on $S^3$ supported by the open book.  Now, since the pages have $(2u+2)$ boundary components, we first attach $(2u+1)$ 1-handles to $D^4$ to get $D^2\times\Sigma_{0,2u+2}$. Then, we attach 2-handles along Legendrian copies of 
boundary parallel curves on $\Sigma_{0,2u+2}$ with framing $\pm 1$, depending on the parity of the Dehn twist. 
The resulting 4-manifold $W$ is given in Figure~\ref{fig:Kirby1}. 

The 1-chain group $C_1(W)$ of $W$ has a basis 
$\{X, Y, Z_1,\dots,Z_{2u-1}\}$
 and $C_2(W)$ has a basis
  $\{ a_1, \dots,a_p,b_1,\dots,b_{p-1},c_1,\dots,c_u,d_1,\dots,d_u \}.$
  The boundary map $D: C_2(W)\to C_1(W)$ is given by 
\begin{align*}
  D(a_j)=X,& & j=1,\ldots p, \\
D(b_j)=Y, & & j=1,\ldots p-1, \\ 
D(c_1)=Z_1-X, D(c_i)=Z_i-Z_{i-1}, & & i=2,\ldots u, \\
D(d_u)=Y-Z_{2u-1}, D(d_i)=Z_{u+i}-Z_{u+i-1}, & & i=1,\ldots u-1.
\end{align*}
Thus, $H_2(W)$ has a basis with generators
$$ \{a_1-a_2,\dots,a_{p-1}-a_p,b_1-b_2,\dots,b_{p-2}-b_{p-1},b_1-\sum_{i=1}^u (c_i+d_i) -a_p\}.$$
Since $\rank H_0=1$, $\rank H_1=0$ and $\rank H_2=2p-2$, we get $\chi(W)=2p-1$.

Note that, $a_j^2=-1=d_j^2$, $b_j^2=1=c_j^2$. So the squares of the basis elements are   $(a_j-a_{j+1})^2=-2$, $(b_j-b_{j+1})^2=2$, $(b_1-\sum_{i=0}^u(c_i+d_i)-a_p)^2=0$. Thus in this basis the intersection matrix is $Q_{\textrm{I}}$ as given in  Appendix~\ref{thematrices}. We also compute  in  Appendix~\ref{theinvariants} that  
$\sigma(W)=\sigma(Q_{\textrm{I}})=0$, and $\det Q_{\textrm{I}}=(-1)^{p-1}$. 

To calculate the square of the first Chern class,  we chose an orientation of the curves and compute the rotation numbers of the curves with respect to the orientation induced from  blackboard.  
Thus we get  $r(a)=0=r(b)$, $r(c_i)=-1$ and $r(d_i)=-1$.  Note that the calculation of $c^2$ is independent of the chosen orientations. Let us denote the cocores of the 2-handles attached along $a_i,b_j,c_k,d_l$ by $A_i,B_j,C_k,D_l$ respectively. 
Then $c(W)$ is Poincar\'e dual to $-(\sum_{i=1}^uC_i+\sum_{j=1}^u D_j)$. This evaluates on the basis above as $w=(0,\dots,2u)^T$. Hence, 
$$c^2(W)=Q_W(PD(c(W)))=w^TQ^{-1}w=\displaystyle\frac{4u^2\cdot(-1)^{p-1}\cdot(p-1)\cdot p}{(-1)^{p-1}}=4u^2p(p-1).$$
Inserting the results of the previous steps in (\ref{eqn:d3invariant})  we get
\begin{equation}
\d3(\xi)=\frac{1}{4}\left(4u^2(p-1)p-2(2p-1)-3\cdot 0\right)+(p+u-1)= u^2p(p-1)+u-\frac{1}{2}.
\end{equation}

\subsection{Overtwisted structures via (II)}
\label{ex:(II)} 
We perform similar calculation for the multilinks (II) given in Figure~\ref{fig:SpliceDiagrams}.
The associated monodromy (\ref{monodromy(II)}) has $q+u$ negative Dehn twists.
After following the same steps to construct the 4-manifold $W$ we find $\chi(W)=q+1$, and as pointed out in Appendix~\ref{theinvariants}, $\sigma(W)=q$. Similarly as before,  we have
    $$ c^2(W)=(2u-1)^2q.$$
Inserting in (\ref{eqn:d3invariant})  we get

\begin{equation}
\displaystyle \d3(\xi)=\frac{1}{4} \left((2u-1)^2q-2(q+1)-3q\right)+(q+u)=u(u-1)q+u-\frac{1}{2}.
\end{equation}

\subsection{Overtwisted structures via (III)}
\label{ex:(III)} 
As for the multilinks (III) in Figure~\ref{fig:SpliceDiagrams},
the associated monodromy (\ref{monodromy(III)}) has $u+1$ negative Dehn twists.
The constructed 4-manifold $W$ has $\chi(W)=5$, and as pointed out in Appendix~\ref{theinvariants}, $\sigma(W)=-2$.  Moreover, we have 
$$c^2(W)=\frac{(2u+1)^2\cdot-6}{-1}=6(2u+1)^2.$$
Inserting in (\ref{eqn:d3invariant})  we get

\begin{equation}
\displaystyle \d3(\xi)=\frac{1}{4}\left(6(2u+1)^2-2\cdot5-3\cdot(-2)\right)+(u+1)=6u(u+1)+u+2-\frac{1}{2}.
\end{equation}
 

\subsection{Overtwisted structures via (I-I)}
\label{ex:(I-I)} 
We consider the graph multilinks (I-I) obtained by splicing two multilinks of type (I), as we have constructed in Figure~\ref{fig:AllSplice}(a).
The monodromy (\ref{monodromy(I-I)}) of the associated open book has $q+u+v-1$ negative Dehn twists.

\begin{figure}[htbp]
	\begin{center}
		\includegraphics[width=1\columnwidth]{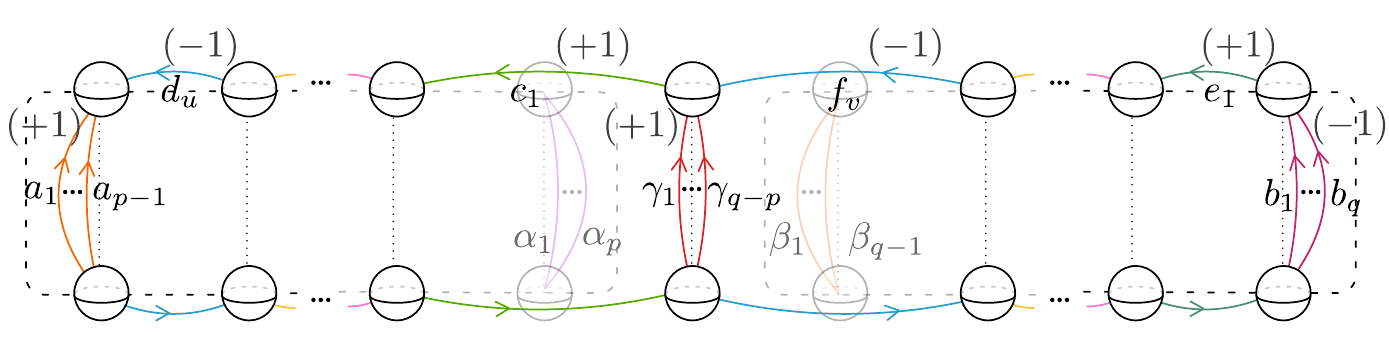}
		\vskip\baselineskip 
		\caption{Kirby diagram for the 4-manifold corresponding to (I-I). The faded ends of the previous diagrams are the  deleted blocks. }
		\label{fig:Kirby11}
	\end{center}
\end{figure}

Since the monodromy is obtained by the monodromies of the splice components, to construct the 4-manifold, we can use the Kirby diagrams for the splice components.  One can see that the Kirby diagram of the spliced multilink can be constructed as follows.  We identify the 1-handles corresponding to the spliced boundary components, thus the 2-handles whose attaching circles corresponds to the Dehn twists along that boundary components cancel.  By means of the new Dehn twist contributions to the monodromy,  we add new 2-handles whose attaching circles are along the identified boundary component.  Consequently,  we see that  the corresponding 4-manifold has the Kirby diagram given in Figure~\ref{fig:Kirby11}.

Furthermore, $H_2(W)$ has a basis with generators
\begin{align*}
& a_1-a_2,\dots,a_{p-2}-a_{p-1}, \gamma_1 - \gamma_2, \dots, \gamma_{k-1}- \gamma_k,  b_1-b_2,\dots,b_{q-1}-b_{q},\\& \hspace{4cm}\gamma_1+(\sum_{i=1}^u c_i+d_i) -a_{p-1}, b_1+(\sum_{i=1}^v e_i+f_i)-\gamma_k.
\end{align*}
Since rank $H_0=1$, rank $H_1=0$ and rank $H_2=2q-2$, we have $\chi(W)=2q-1$.

Note that, $a_j^2=c_j^2=e_j^2=\gamma_j^2=1$ and $b_j^2=d_j^2=f_j^2=-1$. So the squares of the basis elements are $(a_j-a_{j+1})^2=2$, $(\gamma_j-\gamma_{j+1})^2=2$, $(b_j-b_{j+1})^2=-2$, $( \gamma_1+(\sum_{i=1}^u c_i+d_i) -a_{p-1})^2=2$ and $(  b_1+(\sum_{i=1}^v e_i+f_i)-\gamma_k)^2=0$. In this basis the intersection matrix is $Q_{\textrm{I-I}}$ as given in   Appendix~\ref{thematrices}. We compute  in  Appendix~\ref{theinvariants} that  
$\sigma(W)=\sigma(Q_{\textrm{I-I}})=0$, and $\det Q_{\textrm{I-I}}=(-1)^{q-1}$. 

Note that, $r(a)=r(\gamma)=r(b)=0$, $r(c_i)=-1$, $r(d_i)=-1$, $r(e_i)=-1$ and  $r(f_i)=-1$. Therefore, we have, $$c(W)=-\sum_{i=1}^{u}(C_i+D_i)-\sum_{j=1}^v (E_j+F_j).$$ This evaluates on the basis above as $w=(0,\dots,-2u,-2v)^T$. In order to calculate $c^2$, it is sufficient to calculate the inverse of last $2\times2$ block of $Q_{\textrm{I-I}}$. We deduce that $c^2(W)=4u^2p(p-1)+8uvq(p-1)+4v^2q(q-1).$ Explicit calculations can be found in Appendix~\ref{theinvariants}.

Inserting all these results in (\ref{eqn:d3invariant})  we get 
\[\begin{split}
  \displaystyle \d3(\xi)&=\frac{1}{4}(4u^2p(p-1)+8uvq(p-1)+4v^2q(q-1)-2(2q-1)-3\cdot0)+q+u+v-1  \\&=u^2p(p-1)+v^2q(q-1)+2uvq(p-1)+u+v-\frac{1}{2}.  
\end{split}\]

As we have seen, the information about the resulting graph link and its fibration can be deduced from the splice components easily.  In the next example, we will construct a wider family of overtwisted contact structures and observe how the procedure goes on.

\subsection{Overtwisted structures via (I-I-I)}
\label{ex:(I)(I)(I)}

\begin{figure}[htbp]
	\begin{center}
		\includegraphics[width=1\columnwidth]{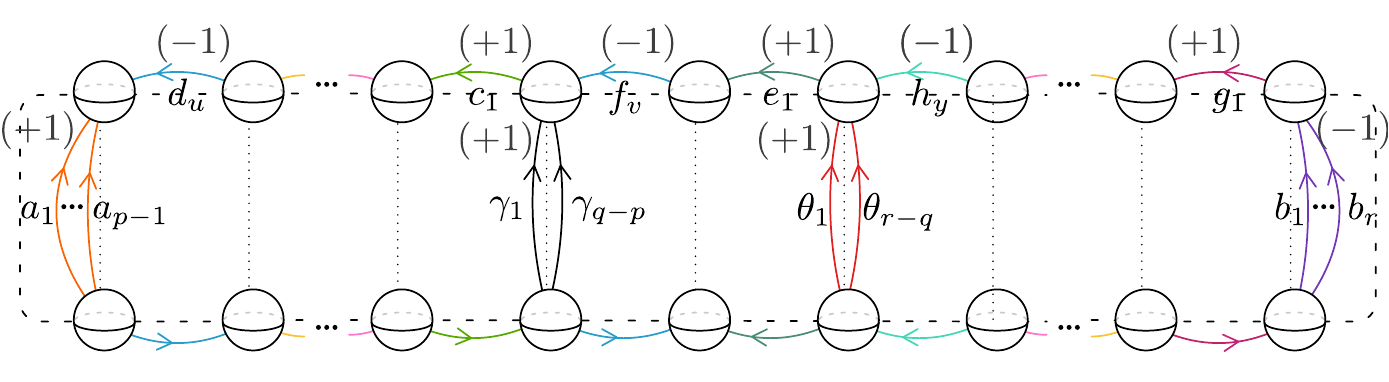}
		\vskip\baselineskip 
		\caption{Kirby diagram for the 4-manifold corresponding to (I-I-I)}
		\label{fig:Kirby111}
	\end{center}
\end{figure}

We consider the graph multilinks (I-I-I) obtained by splicing three multilinks of type (I), as we have constructed in Figure~\ref{fig:Splice111}.
The monodromy  (\ref{monodromy(I-I-I)}) of the associated open book has $r+u+v+w-1$ negative Dehn twists.
By the same arguments as in the previous example, the corresponding 4-manifold has the  Kirby diagram given in Figure~\ref{fig:Kirby111}.


Then $H_2(W)$ has a basis with generators: 
\begin{align*}
&a_1-a_2,\dots,a_{p-2}-a_{p-1}, \gamma_1 - \gamma_2, \dots, \gamma_{q-p-1}- \gamma_{q-p},\theta_1-\theta_2,\dots,\theta_{r-q-1}-\theta_{r-q}, \\
& \hspace{2cm} b_1-b_2,\dots,b_{r-1}-b_{r}, \gamma_1+(\sum_{i=1}^u c_i+d_i) -a_{p-1},\\
&\hspace{2cm} \theta_1+(\sum_{i=1}^v e_i+f_i)-\gamma_{q-p},b_1+(\sum_{i=1}^w g_i+h_i)-\theta_{r-q}.
\end{align*}
Since rank $H_0=1$, rank $H_1=0$ and rank $H_2=2r-2$, we have $\chi(W)=2r-1$.

Note that, $a_j^2=c_j^2=e_j^2=g_j^2=\gamma_j^2=\theta_j^2=1$ and $b_j^2=d_j^2=f_j^2=h_j^2=-1$. So the squares of the basis elements are $(a_j-a_{j+1})^2=2$, $(\gamma_j-\gamma_{j+1})^2=2$, $(\theta_j-\theta_{j+1})^2=2$, $(b_j-b_{j+1})^2=-2$, $( \gamma_1-(\sum_{i=1}^u c_i+d_i) -a_{p-1})^2=2$, $( \theta_1-(\sum_{i=1}^v e_i+f_i) -\gamma_{q-p})^2=2$ and $(  b_1-(\sum_{i=1}^w g_i+h_i)-\theta_{r-q})^2=0$. In this basis the intersection matrix is $Q_{\textrm{I-I-I}}$
as given in Appendix~\ref{thematrices}. We compute  in  Appendix~\ref{theinvariants} that     $\det Q_{\textrm{I-I-I}}=(-1)^{q-1}$. 

As we discussed in the previous example the number of positive eigenvalues is $(p-2)+(q-p-1)+(r-q-1)+3=r-1$ and the number of negative eigenvalues is $(r-1)$. Thus, $\sigma(W)=0$. 
 
 Note that, $r(a)=r(b)=r(\gamma)=r(\theta)=0$, whereas  $r(c_i)=-1$, $r(d_i)=-1$, $r(e_i)=-1$,  $r(f_i)=-1$, $r(g_i)=-1$ and  $r(h_i)=-1$. Therefore, we have $$c(W)=-\sum_{i=1}^{u}(C_i+D_i)-\sum_{j=1}^v (E_j+F_j)-\sum_{j=1}^w (G_j+H_j).$$This evaluates on the basis above as $w=(0,\dots,-2u,-2v, -2w)^T$. In order to calculate $c^2$, it is sufficient to calculate the inverse of the last $3\times3$ block of $Q_{\textrm{I-I-I}}$. The calculations in Appendix~\ref{theinvariants} show 
 $$c^2(W) = 4u^2p(p-1)+4v^2q(q-1)+4w^2r(r-1) +8uvq(p-1)+8uwr(p-1)+8vwr(q-1).$$
Inserting in (\ref{eqn:d3invariant})  we get
\begin{equation}
\begin{split}
 \d3(\xi)&=\frac{1}{4}\left(\begin{split}&4u^2p(p-1)+4v^2q(q-1)+4w^2r(r-1)+8uvq(p-1)\\&+8uwr(p-1)+8vwr(q-1) -2\cdot(2r-1)-3\cdot(0)\end{split}\right)\\&\quad+(r+u+v+w-1)  \\
 &=u^2p(p-1)+v^2q(q-1)+w^2r(r-1)+2uvq(p-1)\\&\quad+2uwr(p-1)+2vwr(q-1)+u+v+w-\frac{1}{2}.  
\end{split}
\label{d3III}
\end{equation}

\subsection{Overtwisted structures via (II-I)}
\label{ex:(I)(II)}

We consider the graph multilinks (II-I) obtained by splicing two multilinks of type (II) (with $q=2$) and (I), as constructed in Figure~\ref{fig:AllSplice}(b).
The monodromy  (\ref{monodromy(II-I)}) of the associated open book has $p+u+v$ negative Dehn twists.
 
 As in the previous examples,  we construct the corresponding 4-manifold and observe that $H_2(W)$ has a basis with generators
 $ a_1-a_2,\dots,a_{p-2}-a_{p-1}, \gamma_1 - \gamma_2, \dots, \gamma_{p-3}- \gamma_{p-2},  b_1-b_2, \gamma_1+(\sum_{i=1}^u c_i+d_i) -a_{p-1}, b_1+\sum_{i=1}^{v-1} e_i+\sum_{i=1}^{v}f_i-\gamma_{p-2}.$ In this basis the intersection matrix is $Q_{\textrm{II-I}}$ as given in the Appendix \ref{thematrices}. There we calculate $\chi(W)=2p-1$,  $\sigma(W)=2$ and det$Q_{\textrm{II-I}}=(-1)^{p}$.  Furthermore, 
 \[\begin{split}
     c^2(W)&=(-2u,-(2v-1))\left(\begin{array}{cc}
p(p-1)&2(p-1)\\
2(p-1)&2
\end{array}\right)(-2u,-(2v-1))^{T}\\
&=4u^2p(p-1)+8v^2+16uv(p-1)-8v-8u(p-1)+2.
\end{split}
\] and we have
\[\begin{split}
  \displaystyle d_3(\xi)&=\frac{1}{4}(4u^2p(p-1)+8v^2+16uv(p-1)-8v-8u(p-1)+2-2(2p-1)-3(2))\\& +p+u+v  \\&=u^2p(p-1)+2v^2+2uv(p-1)+2u(v-1)(p-1)-2v+u+v-\frac{1}{2}.  
\end{split}\]

\subsection{Overtwisted structures via (III-I)}
\label{ex:(I)(III)}

We consider the graph multilinks (III-I) obtained by splicing two multilinks of type (III) and (I), as constructed in Figure~\ref{fig:AllSplice}(c).
 The monodromy  (\ref{monodromy(III-I)}) of the associated open book has $p+u+v-2$ negative Dehn twists.
$H_2(W)$ has a basis with generators $ a_1-a_2,\dots,a_{p-1}-a_{p}, \gamma_1 - \gamma_2, \dots, \gamma_{p-4}- \gamma_{p-3},  b_1-b_2, \gamma_1-(\sum_{i=1}^u c_i+d_i) -a_{p}, b_1-\sum_{i=1}^{v+1} e_i-\sum_{i=1}^{v}f_i-\gamma_{p-3} $.  In this basis the intersection matrix is $Q_{\textrm{III-I}}$ as given in the Appendix \ref{thematrices}.
Similar calculations as before show that $\chi(W)=2p-1$,  $\sigma(W)=-2$,  $\det Q_{\textrm{III-I}}=(-1)^{p}$ and 
\[\begin{split}
     c^2(W)&=(2u,2v+1)\left(\begin{array}{cc}
p(p-1)&2p\\
2p&6
\end{array}\right)(2u,2v+1)^{T}\\
&=4u^2p(p-1)+24v^2+8up(2v+1)+24v+6.
\end{split}
\]

Inserting the results of the previous steps into the formula of $\d3$ invariant, we obtain: 

\[\begin{split}
  \displaystyle d_3(\xi)&=\frac{1}{4}(4u^2p(p-1)+24v^2+8up(2v+1)+24v+6-2(2p-1)-3(-2))\\&+p+u+v-2  \\&=u^2p(p-1)+6v^2+2pu(2v+1)+6v+u+v+2-\frac{1}{2}.  
\end{split}\]

\subsection{Overtwisted structures via (II-III)}
\label{ex:(II)(III)}

Here we consider the graph multilinks (II-III) obtained by splicing two multilinks of type (II)  (with $q=2$) and (III), as constructed in Figure~\ref{fig:AllSplice}(d).
The monodromy  (\ref{monodromy(II-III)}) of the associated open book has $u+v+2$ negative Dehn twists.
Furthermore,  $H_2(W)$ has a basis given as $ a_1-a_2, b_1-b_2, \gamma-\sum_{i=1}^{u-1} c_i-\sum_{i=1}^u d_i -a_{2}, b_1-\sum_{i=1}^{v+1} e_i-\sum_{i=1}^{v}f_i-\gamma$.  In this basis the intersection matrix is $Q_{\textrm{II-III}}$ as given in the Appendix \ref{thematrices}. After similar calculations, we get 
$\chi(W)=5$,  $\sigma(W)=0$,  $\det Q_{\textrm{II-III}}=1$ and 
\[\begin{split}
     c^2(W)&=(2u-1,2v+1) \left(\begin{array}{cc}
2&4\\
4&6\\
\end{array}\right)(2u-1,2v+1)^{T}\\
&=8u^2+24v^2+32uv+8u+8v.
\end{split}
\]
Inserting the results of the previous steps into the formula of $\d3$ invariant,  we obtain: 

\[\begin{split}
  \displaystyle d_3(\xi)&=\frac{1}{4}(8u^2+24v^2+32uv+8u+8v-2\cdot(5)-3\cdot(0))+2+u+v\\&=2u^2+6v^2+4uv+3v+3u-\frac{1}{2}.  
\end{split}\]

\subsection{Proof of the main theorem}
\label{ohfinally}
Finally here we prove our main theorem by showing first that the family of fibered multilinks we obtained by splicing (I-I-I) give us all the overtwisted contact structures with $\d3+\frac{1}{2}>431$ except $\d3+\frac{1}{2}=461$. Then we show that 
all the remaining ones, except for the ones with $\d3 +\frac{1}{2} \in \{4,11,17,19,47,61,79,95,109\}$, are obtained by the other ways of splicing
that we have presented in the previous paragraphs of the present section. We will give a list for that at the end of the section. We do not  know yet if the 9 overtwisted structures 
that we have missed are real algebraic.


Let $d\in\SymZ$ denote the  sum $d_3+1/2$ in Equation~\ref{d3III}:
$$d=u^2p(p-1)+v^2q(q-1)+w^2r(r-1)+2uvq(p-1)+2uwr(p-1)+2vwr(q-1)+u+v+w$$
where the variables are positive integers with the algebraicity condition $p<q<r$. We fix $v=w=1$ once and for all.
We will use the three moves below:
\begin{itemize}
    \item[(i)] Replacing $q$ and $r$ with  $(q+1)$ and $(r-1)$: this increases $d$ by 2.
    \item[(ii)] Replacing $u$  with $(u+2)$ and $r$ with $(r-2)$: this increases $d$ by $4u+12$.
    \item[(iii)] Increasing $r$ by 1: this increases $d$ by $2(r+u+q-1)$.
\end{itemize}

We start from the {\em state} $(p,q,r,u,v,w)=(2,3,r,1,1,1)$. These values give $d=r^2+5r+17$, which is odd.  Any application of the moves above produces an odd number.
First we will tell how to obtain all odd integers greater than 431 (except 461) via these moves.

Starting from the initial state and applying the move (iii) for each $r$ increases the sum by $2r+6$. 
We discuss how to obtain any odd number between $d=r^2+5r+17$ and $d+2r+6=(r+1)^2+5(r+1)+17$ using the first two moves, provided that $r$ is large enough.  

Now starting from the initial state the application of  (ii) $k$ times increases $d$ by $4k^2+12k$ 
Let $k$ be the largest integer satisfying $4k^2+12k<2r+6$.
Note that  we have $k=1$ for $5<r\leq 17$,   $k=2$ for $17<r\leq 33$ and $k=3$ for $33<r\leq 53$. 

Furthermore any odd number between $d+4c^2+12c$ and $d+4(c+1)^2+12(c+1)$ for $0\leq c< k$ can be obtained by applying move (i)
$\frac{8c+16}{2}-1=4c+7$ times. Recall that we have the restriction $q<r$ and that application of moves (i) and (ii) decreases $r$.
Hence in order to obtain all the values in between  we must have $q+4c+7<(r- 2c)-4c-7$, i.e. $r>10c+17$ for any $0\leq c<k$.

When $c=0$,  any odd number between $d$ and $d+16$ can be obtained 
by applying move (i) 7 times. Therefore,  we have the restriction that $q+7<r-7$ hence  $r>17$. 

We have observed above that for $17<r\leq 33$,   the move (ii)  is applied twice.  
Hence for $c=1$ any odd number
between $d+16$ and $d+40$  can be obtained for $r>27$. For  $17<r\leq 27$, there are few values less than $d+40$ that we cannot obtain in this way.  

As for $33<r\leq 53$,  we can apply move (ii) thrice.  Since $r>27$, we have observed above that any sum between $d$ and  $d+40$ can be obtained.  For $c=2$ 
for the numbers between $d+40$ and $d+72$, we must have $r>37$. Hence whenever $37<r\leq 53$ we can obtain any odd number 
between $d$ and $d+72$. For $33<r\leq 37$ we cannot obtain all the numbers in between though.  
For larger $r$ (more precisely for $r>37$) the inequality   $r>10c+17$ is always satisfied  so that we can obtain any odd number
between $d+4c^2+12c$ and $d+4(c+1)^2+12(c+1)$.

Finally in order to obtain any odd number between $d+4k^2+12k$ and $d+2r+6$ via move (i), we must have
$$r-2k-\left(\frac{2r+6-4k^2-12k}{2}-1\right)>3+\frac{2r+6-4k^2-12k}{2}-1,$$ 
i.e. $r<4k^2+10k-7$.  Recall that $k$ is the largest integer satisfying $4k^2+12k<2r+6$.
Comparing these inequalities,  one can  see that  when $r>33$
any odd number in between can be obtained via move (i). 
As a result  we conclude that for $r\geq 18$, i.e. starting from $d=431$  all the odd integers are obtained, except some finitely many missed ones 
for $29\leq r\leq 37$.  Precisely the number of these missed ones is 45. 

Here one can find the exact states that give these missing numbers on a computer. Instead we try to enrich our set of moves in order to obtain most of these 45 numbers.
Indeed,  at the state $(2,3,r-2,3,1,1)$ when we have the sum $d+16$, we increase $p$ and $q$ by 1,  decrease $r$ by 4  to get to 
the state $(3,4,r-6,3,1,1)$ and the sum $d+36$.  Then  applying   the move (i) successively produces the missing numbers between $d+36$ and $d+40$.  Thereby,  we can obtain 18 out of 30 missing odd numbers  between $17<r\leq 27$.  For the 15 missing odd numbers between $29\leq r\leq 37$,  we  replace $p,q,r$ by $p+2$, $q+5$ and $r-11$ at the state $(2,3,r-2,3,1,1)$ to get to the state $(4,8,r-13,3,1,1)$ and the sum $d+62$.  Again,  the application of the move  (i) successively  produces all the odd numbers between $d+62$ and $d+2r+6$. 
For the remaining 12 missing odd numbers  smaller than $d+36$,   at the state $(2,3,r-2,3,1,1)$ we  replace $u$ with $u+4$,  and $r$ with $r-5$.  Application of this move increases the sum by $6u-2r+30$, thus we can obtain all the odd numbers except $d=461$. \\

To obtain even numbers, we start from the state $(2,3,r,2,1,1)$ that gives the even integer $d=r^2+7r+30$.  Then move (iii)  increases the sum by $2r+8$. We will  now obtain any even number between $d$ and $d+2r+8=(r+1)^2+7(r+1)+30$ by applying the first two moves.  
Let  $k$ be the largest integer satisfying $4k^2+16k<2r+8$.
Applying (ii) $k$ times takes us to the state $(2,3,r-2k,2+2k,1,1)$ and increases the value by $4k^2+16k$. 
Each application of (ii), while passing from the step $u+2k$ to $u+2(k+1)$, increases the value by $8k+20$.  
Note that $k=1$ for $6<r\leq 20$,  we have $k=2$ for $20<r\leq 38$ and $k=3$ for $38<r\leq 60$.

Any number between $d+4c^2+16c$ and $d+4(c+1)^2+16(c+1)$ for $0\leq c< k$ can be obtained by applying move (i) $4c+9$ times.  In order to obtain all the sums in between  we must have  $r>10c+21$ for any $0\leq c<k$.
When $c=0$  any even number between $d$ and $d+20$ can be obtained by applying move (i) 9 times for $r>21$. 
We observed above that for $20<r\leq 38$  we apply (ii) twice.  Hence any even number between $d+20$ and $d+48$ (i.e. for $c=1$) can be obtained whenever $r>31$. 
For  $20<r\leq 31$ there are few values less than $d+48$ that we cannot obtain in this way.  

For $38<r\leq 60$  we can apply (ii) thrice.  Any sum between $d$ and $d+48$ can be obtained as discussed in the previous arguments.  For $c=2$, for the numbers between $d+48$ and $d+84$ we must have $r>41$. As a result,  when $41<r\leq 60$, we can obtain any even number between $d$ and $d+84$.  Moreover,  larger $r$ values always satisfy  $r>10c+21$ and we can obtain any even number between $d+4c^2+16c$ and $d+4(c+1)^2+16(c+1)$.

Finally in order to obtain any even number between $d+4k^2+16k$ and $d+2r+8$ via move (i), we must have $r<4k^2+14k-9$.
Checking for the values of $k$,  one can see that the above condition is satisfied for $r>41$.
As a result  we conclude that for $r\geq 17$, i.e. starting from $d=438$  all the even integers are obtained, except some finitely many missed ones. 
Precisely the number of these missed ones is 78.

We have the following additional operations to produce the missed even numbers. To obtain the ones between $d+14$ and $d+20$, 
 at the state $(2,3,r,2,1,1)$ we increase $p$ and $q$ by 1 and decrease $r$ by 3.  The new state $(3,4,r-3,2,1,1)$ gives the sum $d+12$.  Then we apply the move (i) successively to produce all the missing even numbers in 
between.  For the missed even numbers  between $d+48$ and $d+84$,  at the state $(2,3,r-4,6,1,1)$  we decrease $u$ by 4,  $r$ by 14 and increase $p$ and $q$ by 6 to get 
the state $(8, 9,r-18, 2,1,1)$ and the sum $d+72$.  Then applying the move (i) successively produces all the missing even integers between $d+78$ and $d+84$.

Moreover,  at the state $(2,3,r-2,4,1,1)$ with $d+20$,  in order to obtain the missed ones between $d+20$ and $d+48$,  we  decrease $u$ by 2,  $r$ by 7 and increase $p$ and $q$ by 3 which increases the sum by 16.  As before,  we can successively apply the move (i) to obtain the missing ones between $d+36$ and $d+48$.  However,  for the small values of $r$,  18 of the missing even numbers cannot be obtained because of the restriction $q<r$ in each step.  We have realized that 12 of these 18 missing numbers can be produced  by the  application of the move (i) successively  at the states $(2,3,r-7,8,1,1)$ with the sum $d+38$.  Lastly,  one can see that the remaining sums $520, 558, 714, 766, 820, 876$ can be obtained by  the states $(5,6,8,1,2,1)$,  $(4,6,10,1,2,1)$,   $(5,7,10,1,2,1)$,  $(5,7,11,1,2,1)$,  $(5,7,12,1,2,1)$,  $(5,7,13,1,2,1).$

Up to now we have proved that any overtwisted structure with $d_3+1/2 > 431$ (except $461$) can be obtained by (I-I-I) splicing. 
Via computer assistance we find that the ones with all the smaller $d_3$'s (except $4, 11, 17, 19, 47, 61, 79, 95, 109$) 
are obtained via splicings as shown in Table~\ref{tabelo}. This finishes the proof.

\begin{table}[h]
{\tiny
\begin{center} 
\begin{tabular}{ |c||l|}
\hline 
type & $d_3+1/2$ \\
\hline 
\hline 
I	& 3,	7,	10,	13,	21,	26,	31,	43,	50,	55,	57,	73,	82,	91,	111,	122, 133
\\ \hline
II	& 1,	6,	8,	12,	14,	18,	20,	22,	27,	28,	30,	32,	33,	38,	42,	44,	45,	52,	54,	56,	63,	66,	68,	70,	\\
& 75,	84,	86,	93,	102,	104,	124,	134,	156,	182,	189,	208
\\ \hline
III	& 2,	15,	40,	77
\\ \hline
II-I	& 25, 37
\\ \hline
III-I	& 9,	23,	35,	49,	59,	113,	347
\\ \hline
II-III & 5
\\ \hline
& 16,	24,	29,	34,	36,	39,	41,	46,	48,	51,	58,	60,	62,	64,	65,	69,	71,	72,	76,	78,	80,	81,	87,	\\
I-I	& 88, 89,	92,	96,	97,	98,	100,	105,	106,	115,	116,	118,	119,	120,	126,	129,	131,	135,	136,	\\
& 138, 140,	142,	144,	146, 153, 155,	157,	164,	165,	166,	168,	177,	181,	188,	192,	215,	\\
& 246,	249,	256,	275,	313,	358,	387,	461
\\ \hline
I-I-I & all other cases \\
\hline
\end{tabular}
\end{center}
}
\caption{How to obtain the overtwisted structures with $d_3+1/2 \leq 461$ (except $4, 11, 17, 19, 47, 61, 79, 95, 109$). 
Usually there are more than one ways to construct each case. Here we give single samples.}
\label{tabelo}
\end{table}

\bibliography{OTinS3}
\bibliographystyle{amsalpha}

\appendix
\section{Intersection Matrices}
\label{thematrices}
Here,  we  give the intersection matrices of the 4-manifolds that we have constructed in Section~\ref{mainproof}, in the bases we presented there.

Let $J_n$ and $\tilde{J_n}$ be the following matrices:
 $$
 {J_n} =\begin{pmatrix}
2&-1&0&\dots&0\\
-1&2&-1&0&\vdots\\
0&-1&\ddots&\dots&0\\
\vdots&0&\dots&2&-1\\
0&\dots&0&-1&2\\
\end{pmatrix}_{n\times n},
\tilde{J_n} =\begin{pmatrix}
-2&1&0&\dots&0\\
1&-2&1&0&\vdots\\
0&1&\ddots&\dots&0\\
\vdots&0&\dots&-2&1\\
0&\dots&0&1&-2\\
\end{pmatrix}_{n\times n}
$$ 

Then the intersection matrices for (I), (II) and (III) are as follows:
{\small
$$
Q_{\textrm{I}}=
\left(
\begin{array}{ccc|ccc|c}
  &&&&&& \\ 
  &\tilde{J}_{p-1}&&&&&\\
  &&&&&&-1\\
  \hline
  &&&&&&1\\
  &&&& J_{p-2}&&\\
  &&&&&&\\
  \hline
  &&-1&1&&&0\\
 \end{array}
 \right),
 Q_{\textrm{II}}=
\left(
\begin{array}{ccc|c}
  &&& \\ 
  &J_{q-1}&&\\
  &&&1\\
  \hline
  &&1&1\\
 \end{array}
 \right),
Q_{\textrm{III}}=\begin{pmatrix}
-2&1&0&0\\
1&-2&0&-1\\
0&0&-2&-1\\
0&-1&-1&-1\\
\end{pmatrix}.
$$
}

The intersection matrices for the 4-manifolds obtained for the spliced graph multilinks  (I-I) and (I-I-I) are respectively as follows.
Here $a=-1$ if $q=p+1$ and $a=0$ for $q>p+1$; $b=-1$ if $r=q+1$ and $b=0$ for $r>q+1$.

{\small
\[
Q_{\textrm{I-I}}=
\left (
\begin{array}{ccc|ccc|ccc|cc}
  &&&&&&&&&& \\ 
  &J_{p-2}&&&&&&&&&\\
  &&&&&&&&&1&\\
  \hline
  &&&&&&&&&1&\\
  &&&&J_{q-p-1}&&&&&&\\
  &&&&&&&&&&1\\
  \hline
  &&&&&&&&&&-1\\
  &&&&&&&\tilde{J}_{q-1}&&&\\
  &&&&&&&&&&\\
  \hline
 &&1&1&&&&&&2&a\\
 &&&&&1&-1&&&a&0\\
 \end{array}
 \right),
\]

$$
Q_{\textrm{I-I-I}}=
\left (
\begin{array}{ccc|ccc|ccc|ccc|ccc}
  &&&&&&&&&&&&&& \\ 
  &J_{p-2}&&&&&&&&&&&&&\\
  &&&&&&&&&&&&1&&\\
  \hline
  &&&&&&&&&&&&1&&\\
  &&&&J_{q-p-1}&&&&&&&&\\
  &&&&&&&&&&&&&1&\\
  \hline
  &&&&&&&&&&&&&1&\\
  &&&&&&&J_{r-q-1}&&&&&&&\\
  &&&&&&&&&&&&&&1\\
   \hline
  &&&&&&&&&&&&&&-1\\
  &&&&&&&&&&\tilde{J}_{r-1}&&&&\\
  &&&&&&&&&&&&&\\
  \hline
 &&1&1&&&&&&&&&2&&a\\
 &&&&&1&1&&&&&&&2&b\\
  &&&&&&&&1&-1&&&a&b&0\\
 \end{array}
 \right).
$$
} 

Finally here are the intersection matrices for the 4-manifolds obtained for the spliced graph multilinks  (II-I) and (III-I) respectively.  In the former matrix, $a=1$ if $p=3$ and $a=0$ for $p\geq 3$. In the latter matrix  $a=-1$ if $p=4$ and $a=0$ for $p\geq 4$.
{\small
$$
Q_{\textrm{II-I}}=
\left(
\begin{array}{ccc|ccc|c|cc}
  &&&&&&&&\\ 
  &J_{p-2}&&&&&&&\\
  &&&&&&&1&\\
  \hline
  &&&&&&&-1&\\
  &&&&\tilde{J}_{p-3}&&&&\\
  &&&&&&&&-1\\
  \hline
  &&&&&&2&&1\\
  \hline
  &&1&-1&&&&0&a\\
  &&&&&-1&1&a&-1\\
 \end{array}
 \right),
$$
  
$$Q_{\textrm{III-I}}=
\left(
\begin{array}{ccc|ccc|c|cc}
  &&&&&&&&\\ 
  &\tilde{J}_{p-1}&&&&&&&\\
  &&&&&&&-1&\\
  \hline
  &&&&&&&1&\\
  &&&&J_{p-4}&&&&\\
  &&&&&&&&1\\
  \hline
  &&&&&&-2&&-1\\
  \hline
  &&-1&1&&&&0&a\\
  &&&&&1&-1&a&1\\
 \end{array}
 \right).
$$
} 
\section{Determinant, signature and $c^2$ computation}
\label{theinvariants}
Here, we give detailed calculations for the results about the intersection matrices we used in Section~\ref{mainproof}.
For practical reference we summarize these results in Table~\ref{tab.sign.det}.
\begin{table}[h]
{\tiny
\begin{center} 
\begin{tabular}{ |c||l||l|l|}
\hline 
 & $\sigma$ & $\det$ & $c^2$\\
\hline 
\hline 
$Q_{\textrm{I}}$ & 0 & $(-1)^{p-1}$ & $4u^2 p(p - 1)$
\\ \hline
$Q_{\textrm{II}}$ & $q$ & 1 & $(2u - 1)^2 q$
\\ \hline
$Q_{\textrm{III}}$ & $-2$ & $-1$ & $6(2u + 1)^2$
 \\ \hline
$Q_{\textrm{I-I}}$ & $0$ & $(-1)^{q-1}$ & $4u^2 p(p - 1) + 8uvq(p - 1)$
 \\ &&& $ + 4v^2 q(q - 1)$
 \\ \hline
$Q_{\textrm{I-I-I}}$ & $0$ & $(-1)^{q-1}$ & $4u^2 p(p-1)+4v^2 q(q-1)+4w^2 r(r-1)$
 \\ &&& $ +8uvq(p-1)+8uwr(p-1)+8vwr(q-1)$
 \\ \hline
$Q_{\textrm{II-I}}$ & $2$ & $(-1)^p$ & $4u^2 p(p - 1) + 8v^2 + 16uv(p - 1) $
 \\ &&& $ - 8v - 8u(p - 1) + 2$
 \\ \hline
$Q_{\textrm{III-I}}$ & $0$ & $1$ & $8u^2 + 24v^2 + 32uv + 8u + 8v$
 \\ \hline
$Q_{\textrm{II-III}}$ & $0$ & $1$ & $8u^2 + 24v^2 + 32uv + 8u + 8v$
 \\
\hline
\end{tabular}
\end{center}
}
\caption{Signatures and determinants of the intersection matrices and the corresponding $c^2$}
\label{tab.sign.det}
\end{table}

\subsection{The matrix $Q_{\textrm{I}}$}
First we compute the diagonalization of the intersection matrix $Q_{\textrm{I}}$ above. Consider the lower triangular matrix $S_n$ with its $ij$ entry ($i\geq j$)
being equal to $j/i$.
It can be easily seen that $J_n=S_n D_n S_n^T$ and $\tilde{J}_n=S_n(-D_n)S_n^T$ where $D_n=\diag(2,\frac{3}{2}, \cdots, \frac{n+1}{n})$.  
It follows that $\det J_n =n+1 \text{ and }\det \tilde{J}_n = (-1)^n(n+1).$

Thus  we see that $Q_{\textrm{I}}=SDS^T$ where $D=\diag(-2,-\frac{3}{2}, \cdots, -\frac{p}{p-1},2,\frac{3}{2},\cdots, \frac{p-1}{p-2},\frac{1}{p(p-1)})$ and 
{\small
\[S=
\left(
\begin{array}{ccc|ccc|c}
  &&&&&& \\ 
  &S_{p-1}&&&&&\\
  &&&&&&-1\\
  \hline
  &&&&&&1\\
  &&&& S_{p-2}&&\\
  &&&&&&\\
  \hline
  -\frac{1}{p}&\cdots&-\frac{p-1}{p}&-\frac{p-2}{p-1}&\cdots&-\frac{1}{p-1}&1\\
 \end{array}
 \right)
 \]
} 
We conclude that $\sigma(Q_{\textrm{I}})=0$ and $\det Q_{\textrm{I}} =(-1)^{p-1}$.  

Similarly, signatures and determinants of the intersection matrices for (II) and (III) can be calculated as given in Table~\ref{tab.sign.det}. 

\subsection{The matrix $Q_{\textrm{I-I}}$}
As for the intersection matrix $Q_{\textrm{I-I}}$ above for the splicing (I-I), we first show that  $\det Q_{\textrm{I-I}}=(-1)^{q-1}$.  
Then we compute the signature and $c^2$.

Let $J_n^{\prime}$ (respectively $J_n^{\prime\prime}$) be the matrix obtained by removing the last (respectively the first) column of the matrix $J_n$.

We assume that $q>p+1$, i.e. $a=0$ in  $Q_{\textrm{I-I}}$; it can be shown that the results are the same in the case  $q=p+1$. 
Now we calculate $\det Q_{\textrm{I-I}}$ via its last row.
We observe that it is equal to $(-1)^{q-1}$ times the sum of the  determinants below  
{\small
\[
\begin{gathered}
{ \left| \;
\begin{array}{cc|cc|cc|cc}
  &&&&&&& \\ 
  J_{p-2}&&&&&&&\\
  &&&&&&1&\\
  \hline
  &&&&&&1&\\
  &&&J_{q-p-1}^{\prime}&&&&\\
  &&&&&&&1\\
  \hline
  &&&&&&&-1\\
  &&&&&\tilde{J}_{q-1}&&\\
  &&&&&&&\\
  \hline
 &1&1&&&&2&\\
 \end{array}
 \; \right|
 }
 +{\left| \; 
\begin{array}{cc|cc|cc|cc}
  &&&&&&& \\ 
  J_{p-2}&&&&&&&\\
  &&&&&&1&\\
  \hline
  &&&&&&1&\\
  &&&J_{q-p-1}&&&&\\
  &&&&&&&1\\
  \hline
  &&&&&&&-1\\
  &&&&&\tilde{J}_{q-1}^{\prime\prime}&&\\
  &&&&&&&\\
  \hline
 &1&1&&&&2&\\
 \end{array}
  \; \right|
}
\end{gathered}\]
} 

Now we move the last column of each matrix above to the positions of the removed columns, i.e. in the first matrix we move the $(2q-3)^{rd}$ column to the  $(q-3)^{rd}$  position and in the second matrix to the $(q-2)^{nd}$  position. These row exchanges multiply the determinants by $(-1)^{(2q-3-q+3)}$ and  $(-1)^{(2q-3-q+2)}$, respectively.  Since 
\[\left| \;
\begin{array}{cc|c}
&J_{n}^{\prime}&\\ 
\hline
&&1\\
\end{array}
\right| = \det J_{n-1}, \quad
\left| \;
\begin{array}{c|cc}
-1&&\\
\hline
&& \\
&&\tilde{J}_{n}^{\prime\prime}\\
\end{array}
\right| =- \det \tilde{J}_{n-1}=(-1)^nn,
\]
we have 
\begin{equation*}
    \begin{split}
    \det Q_{\textrm{I-I}}&=\quad(-1)^{q-1}\left(
    \begin{split}
  &\quad(-1)^{q-1}\cdot\det J_{p-2}^{\prime}\cdot \det J_{q-p-1}^{\prime}\cdot \tilde{J}_{q-1} \\
  &+(-1)^{q-1}\cdot\det J_{p-2}\cdot \det J_{q-p-1}^{\prime\prime\prime}\cdot\det \tilde{J}_{q-1}\\
  &+(-1)^q\cdot 2\cdot\det J_{p-2}\cdot \det J_{q-p-1}^{\prime}\cdot \det \tilde J_{q-1}
  \end{split} \right)\\
  &\quad+(-1)^{q-1}\left(
    \begin{split}
  &\quad(-1)^{q}\cdot\det J_{p-2}^{\prime}\cdot \det J_{q-p-1}\cdot \tilde{J}_{q-1}^{\prime\prime} \\
  &+(-1)^{q}\cdot\det J_{p-2}\cdot \det J_{q-p-1}^{\prime\prime}\cdot\det \tilde{J}_{q-1}^{\prime\prime}\\
  &+(-1)^{q-1}\cdot 2\cdot\det J_{p-2}\cdot \det J_{q-p-1}\cdot \det \tilde J_{q-1}^{\prime\prime}
  \end{split} \right)\\
  &=(-1)^qq(q-2)-(-1)^q(q-1)^2\\
  &=(-1)^{q-1}.
  \end{split} 
  \end{equation*}

 Now we compute the signature of $Q_{\textrm{I-I}}$. This matrix is in the form
  \[
\left[
\begin{array}{c|c}
 A& B^{T} \\
\hline
B& C
\end{array}
\right]
\] 
where $A$, $B$, $C$ are $(2q-4)\times(2q-4)$, $2\times(2q-4)$ and  $2\times2$ symmetric matrices respectively.
Let $S_1$ and $S_2$ be the orthogonal matrices that diagonalize $A$ and $C-BA^{-1}B^{T}$ respectively. Define 
   \[S=
\left[
\begin{array}{c|c}
 S_1& 0 \\
\hline
-S_2BA^{-1}& S_2
\end{array}
\right]
.\]  It can be  seen easily that 
 \[SQS^{T}=
\left[
\begin{array}{c|c}
 S_1AS_1^{T}& 0 \\
\hline
0& S_2(C-BA^{-1}B^{T})S_2^{T}
\end{array}
\right]
.\] Hence 
$$\sigma(Q_{\textrm{I-I}})=\sigma(SQ_{\textrm{I-I}}S^{T})=\sigma( S_1AS_1^{T})+\sigma(S_2(C-BA^{-1}B^{T})S_2^{T})=\sigma(A)+\sigma(C-BA^{-1}B^{T}).$$
We know that $J_n$ and $\tilde{J}_n$ are diagonalizable and are positive definite and negative definite  respectively. Therefore, $A$ has $(p-2)+(q-p-1)=q-3$ positive and $q-1$ negative eigenvalues. Moreover it is easy to observe that $C-BA^{-1}B^{T}$ is positive definite, hence has 2 positive eigenvalues. Hence  $\sigma(Q_{\textrm{I-I}})=0$.

Now we compute $c^2$.  As we have observed,  the basis of $H_2(W;\mathbb{Z})$ given in Section~\ref{ex:(I)}, $c(W)$ evaluates on  as $w=(0,\dots,-2u,-2v)^T$. Hence, in order to calculate $c^2$, it is sufficient to calculate the inverse of the last $2\times2$ block of $Q_{\textrm{I-I}}$. Let $D$  denote this matrix and $d_{ij}$ be its $(i,j)$ entry. We claim
$$
d_{11}=\frac{\cofac_{11}}{\det Q_{\textrm{I-I}}}=p(p-1),\quad d_{12}=d_{21}=q(p-1) \mbox{ and }  d_{22}=q(q-1). 
$$
To prove these we compute the cofactors explicitly. First,
\begin{align*}
\cofac_{11}&=
\left |
\begin{array}{cc|ccc|cc|c}
  &&&&&&& \\ 
  J_{p-2}&&&&&&&\\
  &&&&&&&\\
  \hline
  &&&&&&&\\
  &&&J_{q-p-1}&&&&\\
  &&&&&&&1\\
  \hline
  &&&&&&&-1\\
  &&&&&&\tilde{J}_{q-1}&\\
  &&&&&&&\\
  \hline
 &&&&1&-1&&0\\
 \end{array}
 \right |\\
&=(-1)^q(-1)^{q-1}\det J_{p-2}\cdot\det J_{q-p-1}^{\prime}\cdot\det\tilde{J}_{q-1}\\ 
&\quad +(-1)^q(-1)^{q-2}\det J_{p-2}\cdot\det J_{q-p-1}\cdot\det\tilde{J}_{q-1}^{\prime\prime}\\
&=(-1)^qq(p-1)(q-p-1)+(-1)^{q-1}(q-1)(p-1)(q-p)\\
&=(-1)^{q-1}p(p-1).
\end{align*}
The following determinants are evaluated by induction:
\[\overline{J_{n}}= \left| \;
\begin{array}{c|cc}
&&J_{n}^{\prime\prime}\\
\hline
1&&
\end{array}
\right| =(-1)^2\overline{J_{n-1}}=1, \qquad \textrm{ and }
\left| \;
\begin{array}{c|cc}
&&J_{n}\\
\hline
1&&
\end{array}
\right|=0. \] 
Therefore,
\begin{align*}
&\cofac_{12}=  \\
 &(-1)^{p-1} \left| \;
\begin{array}{cc|cc|cc|c}
  &&&&&& \\ 
  J_{p-2}^{\prime}&&&&&&\\
  &&&&&&\\
  \hline
  &&&&&&\\
  &&&J_{q-p-1}&&&\\
  &&&&&&1\\
  \hline
  &&&&&&-1\\
  &&&&&\tilde{J}_{q-1}&\\
  &&&&&&\\
 \end{array}
 \; \right| \quad + (-1)^{p}\left| \;
\begin{array}{cc|cc|cc|c}
  &&&&&& \\ 
  J_{p-2}&&&&&&\\
  &&&&&&\\
  \hline
  &&&&&&\\
  &&&J_{q-p-1}^{\prime\prime}&&&\\
  &&&&&&1\\
  \hline
  &&&&&&-1\\
  &&&&&\tilde{J}_{q-1}&\\
  &&&&&&\\
 \end{array}
\; \right|
\\
&\quad=0+(-1)^p(-1)^{p-1}(-1)^{q-1}q(p-1)\\
&\quad=(-1)^qq(p-1).
\end{align*}
\begin{align*}
\cofac_{22}&=
\left |
\begin{array}{cc|ccc|cc|c}
  &&&&&&& \\ 
  J_{p-2}&&&&&&&\\
  &&&&&&&1\\
  \hline
  &&&&&&&1\\
  &&&J_{q-p-1}&&&&\\
  &&&&&&&\\
  \hline
  &&&&&&&\\
  &&&&&&\tilde{J}_{q-1}&\\
  &&&&&&&\\
  \hline
 &1&1&&&&&2\\
 \end{array}
 \right | \\
&=(-1)^{p-1}(-1)^{p}\det J_{p-2}^{\prime}\cdot\det J_{q-p-1}\cdot\det\tilde{J}_{q-1}\\&\quad+(-1)^p(-1)^{p-1}\det J_{p-2}\cdot\det J_{q-p-1}^{\prime\prime}\cdot\det\tilde{J}_{q-1}\\&\quad+2\det J_{p-2}\cdot\det J_{q-p-1}\cdot\det\tilde{J}_{q-1}\\
&=(-1)^{q-1}q(q-1).
\end{align*}
As a result, we conclude that
\begin{align*}
c^2(W)&=(-2u,-2v)\left(\begin{array}{cc}
p(p-1)&q(p-1)\\
q(p-1)&q(q-1)
\end{array}\right)(-2u,-2v)^{T}\\
&=4u^2p(p-1)+8uvq(p-1)+4v^2q(q-1).
\end{align*}

\subsection{The matrix $Q_{\textrm{I-I-I}}$}
Now we consider  the splicing (I-I-I) and the associated  intersection matrix $Q_{\textrm{I-I-I}}$ given in Appendix~\ref{thematrices}.  We omit the calculations of the cases when $a,b$ are nonzero since they give the same results.  
Calculating the determinant of the intersection matrix with $a=b=0$ with respect to its last three rows as in the previous example we have
\begin{align*}
\det Q_{\textrm{I-I-I}}
&=(-1)^{r-1}r(r-q-1)(q-2)+(-1)^{r-1}r(r-q-2)(q-1) \\
&\hspace{1cm}-2(-1)^{r-1}r(r-q-1)(q-1)\\
&+(-1)^r(r-1)(r-q)(q-2)+(-1)^r(r-1)(r-q-1)(q-1)\\
&\hspace{1cm}-2(-1)^r(r-1)(r-q)(q-1)\\
 &=(-1)^{r-1}.
 \end{align*}
Moreover  $c(W)$ evaluates on the given basis of $H_2(W;\mathbb{Z})$ as $w=(0,\dots,-2u,-2v, -2y)^T$. In order to calculate $c^2$, it is sufficient to calculate the inverse of the last $3\times3$ block $D$ of $Q_{\textrm{I-I-I}}$. We have
\[
D= \left(\begin{array}{ccc}
p(p-1)&q(p-1)&r(p-1)\\
q(p-1)&q(q-1)&r(q-1)\\
r(p-1)&r(q-1)&r(r-1)\\
\end{array}\right)
.\]
We see that
{\small
\begin{equation*}
\cofac_{11}= 
\left| \;
\begin{array}{cc|cc|cc|cc|cc}
  &&&&&&&&& \\ 
  J_{p-2}&&&&&&&&&\\
  &&&&&&&&&\\
  \hline
  &&&&&&&&&\\
  &&J_{q-p-1}&&&&&&&\\
  &&&&&&&&1&\\
  \hline
  &&&&&&&&1&\\
  &&&&J_{r-q-1}&&&&\\
  &&&&&&&&&1\\
   \hline
  &&&&&&&&&-1\\
  &&&&&&&\tilde{J}_{r-1}&&\\
  &&&&&&&&&\\
  \hline
  &&&1&1&&&&2&\\
  &&&&&1&-1&&&0\\
 \end{array}
 \; \right|
\end{equation*}
} 
which can be calculated with respect to the last two rows yielding
\begin{align*}
\cofac_{11}
&=\det J_{p-2}\cdot\det J_{q-p-1}^{\prime}\cdot\det J_{r-q-1}^{\prime}\cdot\det \tilde{J}_{r-1}\\ 
&\hspace{1cm}+\det J_{p-2}\cdot\det J_{q-p-1}\cdot\det J_{r-q-1}^{\prime\prime\prime}\cdot\det \tilde{J}_{r-1} \\
&\hspace{2cm}-2\det J_{p-2}\cdot\det J_{q-p-1}\cdot\det J_{r-q-1}^{\prime}\cdot\det \tilde{J}_{r-1}\\
&+\det J_{p-2} \cdot \det J_{q-p-1}^{\prime} \cdot \det J_{r-q-1} \cdot \det \tilde{J}_{r-1}^{\prime\prime}\\
&\hspace{1cm}+\det J_{p-2} \cdot \det J_{q-p-1} \cdot \det J_{r-q-1}^{\prime\prime} \cdot \det \tilde{J}_{r-1}^{\prime\prime}\\
&\hspace{2cm}-2\det J_{p-2} \cdot \det J_{q-p-1}\cdot \det J_{r-q-1} \cdot \det \tilde{J}_{r-1}^{\prime\prime}\\
 &=(-1)^{r-1}r(p-1)(p-r+1)+(-1)^r(r-1)(p-q)(p-r)\\
 &=(-1)^{r-1}p(p-1).
\end{align*}
Meanwhile,
{\small
\begin{equation*}
\cofac_{22}=\left| \;
\begin{array}{cc|cc|cc|cc|cc}
  &&&&&&&&& \\ 
  J_{p-2}&&&&&&&&&\\
  &&&&&&&&1&\\
  \hline
  &&&&&&&&1&\\
  &&J_{q-p-1}&&&&&&&\\
  &&&&&&&&&\\
  \hline
  &&&&&&&&&\\
  &&&&J_{r-q-1}&&&&\\
  &&&&&&&&&1\\
   \hline
  &&&&&&&&&-1\\
  &&&&&&&\tilde{J}_{r-1}&&\\
  &&&&&&&&&\\
  \hline
  &1&1&&&&&&2&\\
  &&&&&1&-1&&&0\\
 \end{array}
  \; \right|
\end{equation*}
} 
which can be calculated with respect to the last two rows so that
\begin{align*}
\cofac_{22}
&=\det J_{p-2}^{\prime}\cdot\det J_{q-p-1}\cdot\det J_{r-q-1}^{\prime}\cdot\det \tilde{J}_{r-1}\\ 
&\hspace{1cm}+\det J_{p-2}\cdot\det J_{q-p-1}^{\prime\prime}\cdot\det J_{r-q-1}^{\prime}\cdot\det \tilde{J}_{r-1} \\
&\hspace{2cm}-2\det J_{p-2}\cdot\det J_{q-p-1}\cdot\det J_{r-q-1}^{\prime}\cdot\det \tilde{J}_{r-1} \\
&+\det J_{p-2}'\cdot \det J_{q-p-1}\cdot \det J_{r-q-1} \cdot \det \tilde{J}_{r-1}^{\prime\prime}\\
&\hspace{1cm}+\det J_{p-2} \cdot \det J_{q-p-1}^{\prime\prime} \cdot \det J_{r-q-1}\cdot \det \tilde{J}_{r-1}^{\prime\prime}\\
&\hspace{2cm}-2\det J_{p-2} \cdot \det J_{q-p-1}\cdot \det J_{r-q-1} \cdot \det \tilde{J}_{r-1}^{\prime\prime}\\
 &=(-1)^{r}q((p-2)(q-p)+(p-1)(q-p-1)-2(p-1)(q-p))\\
 &=(-1)^{r-1}q(q-1).
\end{align*}
Similarly,
{\small
\begin{equation*}
\cofac_{33}=\left| \;
\begin{array}{cc|cc|cc|cc|cc}
  &&&&&&&&& \\ 
  J_{p-2}&&&&&&&&&\\
  &&&&&&&&1&\\
  \hline
  &&&&&&&&1&\\
  &&J_{q-p-1}&&&&&&&\\
  &&&&&&&&&1\\
  \hline
  &&&&&&&&&1\\
  &&&&J_{r-q-1}&&&&\\
  &&&&&&&&&\\
   \hline
  &&&&&&&&&\\
  &&&&&&&\tilde{J}_{r-1}&&\\
  &&&&&&&&&\\
  \hline
  &1&1&&&&&&2&\\
  &&&1&1&&&&&2\\
 \end{array}
  \; \right|
\end{equation*}
} 
which can be calculated with respect to the last two rows so that
\begin{align*}
\cofac_{33}
&=\det J_{p-2}^{\prime}\cdot\det J_{q-p-1}^{\prime}\cdot\det J_{r-q-1}\cdot\det \tilde{J}_{r-1}\\ 
&\hspace{1cm}+\det J_{p-2}\cdot\det J_{q-p-1}^{\prime\prime\prime}\cdot\det J_{r-q-1}\cdot\det \tilde{J}_{r-1} \\
&\hspace{2cm}-2\det J_{p-2}\cdot\det J_{q-p-1}^{\prime}\cdot\det J_{r-q-1}\cdot\det \tilde{J}_{r-1} \\
&+\det J_{p-2}^{\prime} \cdot \det J_{q-p-1} \cdot \det J_{r-q-1}^{\prime\prime} \cdot \det \tilde{J}_{r-1}\\
&\hspace{1cm}+\det J_{p-2} \cdot \det J_{q-p-1}^{\prime\prime}  \cdot \det J_{r-q-1}^{\prime\prime} \cdot \det \tilde{J}_{r-1}\\
&\hspace{2cm}-2\det J_{p-2} \cdot \det J_{q-p-1}\cdot \det J_{r-q-1}^{\prime\prime}  \cdot \det \tilde{J}_{r-1}\\
&+2\det J_{p-2}^{\prime}\cdot \det J_{q-p-1} \cdot \det J_{r-q-1}\cdot \det \tilde{J}_{r-1}\\ 
&\hspace{1cm}+2\det J_{p-2} \cdot \det J_{q-p-1}^{\prime\prime}  \cdot \det J_{r-q-1} \cdot \det \tilde{J}_{r-1}\\
&\hspace{2cm}+4\det J_{p-2} \cdot \det J_{q-p-1}\cdot \det J_{r-q-1}\cdot \det \tilde{J}_{r-1} \\
 &=(-1)^{r-1}r(r-1).
\end{align*}
Note also that
\[ \overline{J_m} =\left| \;
\begin{array}{c|cc}
1&&\\
\hline
&&J_{m}^{\prime}
\end{array}
 \; \right| =(-1)^{m-1},
\] which follows by induction. Thus
{\small
\begin{align*}
\cofac_{12}&=(-1)^{r-1}(-1)^{q-1} \left| \;
\begin{array}{cc|cc|cc|cc|cc}
  &&&&&&&&& \\ 
  J_{p-2}&&&&&&&&&\\
  &&&&&&&&1&\\
  \hline
  &&&&&&&&1&\\
  &&J_{q-p-1}^{\prime}&&&&&&&\\
  &&&&&&&&&\\
  \hline
  &&&&&&&&&\\
  &&&&&J_{r-q-1}^{\prime}&&&&\\
  &&&&&&&&&1\\
   \hline
  &&&&&&&&&-1\\
  &&&&&&&\tilde{J}_{r-1}&&\\
  &&&&&&&&&\\
 \end{array}
  \; \right|\\
 & \quad\text{  }+(-1)^{r-1}(-1)^{q}\left| \;
\begin{array}{cc|cc|cc|cc|cc}
  &&&&&&&&& \\ 
  J_{p-2}&&&&&&&&&\\
  &&&&&&&&1&\\
  \hline
  &&&&&&&&1&\\
  &&&J_{q-p-1}&&&&&&\\
  &&&&&&&&&\\
   \hline
  &&&&&&&&&\\
  &&&&&J_{r-q-1}^{\prime\prime\prime}&&&&\\
  &&&&&&&&&1\\
  \hline
  &&&&&&&&&-1\\
  &&&&&&&\tilde{J}_{r-1}&&\\
  &&&&&&&&&\\
 \end{array}
  \; \right|
\end{align*}

\begin{align*}
&\quad\quad\quad+(-1)^{r-1}(-1)^{q-1} \left| \;
\begin{array}{cc|cc|cc|cc|cc}
  &&&&&&&&& \\ 
  J_{p-2}&&&&&&&&&\\
  &&&&&&&&1&\\
  \hline
  &&&&&&&&1&\\
  &&J_{q-p-1}^{\prime}&&&&&\\
  &&&&&&&&&\\
  \hline
  &&&&&&&&&\\
  &&&&J_{r-q-1}&&&\\
  &&&&&&&&&1\\
  \hline
  &&&&&&&&&-1\\
  &&&&&&&\tilde{J}_{r-1}^{\prime\prime}&&\\
  &&&&&&&&&\\
 \end{array}
  \; \right|\\
&\quad\quad\quad\quad+(-1)^{r-1}(-1)^q\left| \;
\begin{array}{cc|cc|cc|cc|cc}
  &&&&&&&&& \\ 
  J_{p-2}&&&&&&&&&\\
  &&&&&&&&1&\\
  \hline
  &&&&&&&&1&\\
  &&&J_{q-p-1}&&&&&&\\
  &&&&&&&&&\\
    \hline
  &&&&&&&&&\\
  &&&&&J_{q-r-1}^{\prime\prime}&&&&\\
  &&&&&&&&&1\\
  \hline
  &&&&&&&&&-1\\
  &&&&&&&\tilde{J}_{r-1}^{\prime\prime}&&\\
  &&&&&&&&&\\
 \end{array}
  \; \right|
\\
&\quad\quad=(-1)^{r-1}r(p-1)(r-q-1)+0+(-1)^r(r-1)(p-1)(r-q)+0\\
&\quad\quad=-(-1)^{r-1}q(p-1).
\end{align*}
} 
Similarly,
{\small
\begin{equation*}
\cofac_{13}=\left| \;
\begin{array}{cc|cc|cc|cc|cc}
  &&&&&&&&& \\ 
  J_{p-2}&&&&&&&&&\\
  &&&&&&&&1&\\
  \hline
  &&&&&&&&1&\\
  &&J_{q-p-1}&&&&&&&\\
  &&&&&&&&&1\\
  \hline
  &&&&&&&&&1\\
  &&&&J_{r-q-1}&&&&\\
  &&&&&&&&&\\
   \hline
  &&&&&&&&&\\
  &&&&&&&\tilde{J}_{r-1}&&\\
  &&&&&&&&&\\
  \hline
  &1&1&&&&&&2&\\
  &&&&&1&-1&&&2\\
 \end{array}
  \; \right|
\end{equation*}
} 
which can be calculated with respect to the last two rows yielding
\begin{align*}
\cofac_{13}
&=-(-1)^{r-1}(-1)^{p-1}\det J_{p-2}\cdot\det \overline{J_{q-p-1}^{\prime}}\cdot\det \overline{J_{r-q-1}^{\prime}}\cdot\det \tilde{J}_{r-1}\\ 
&\hspace{1cm}+(-1)^{r-1}(-1)^{p-1}\det J_{p-2}\cdot\det \overline{J_{q-p-1}}\cdot\det \overline{J_{r-q-1}^{\prime\prime\prime}}\cdot\det \tilde{J}_{r-1}\\
&\hspace{2cm}-(-1)^{r-1}(-1)^q2\det J_{p-2}\cdot\det J_{q-p-1}\cdot\det \overline{J_{r-q-1}^{\prime}}\cdot\det \tilde{J}_{r-1}  \\
&-(-1)^{r-1}(-1)^{p-1}\det J_{p-2} \cdot \det \overline{J_{q-p-1}'} \cdot \det \overline{J_{r-q-1}} \cdot \det \tilde{J}_{r-1}^{\prime\prime}\\
&\hspace{1cm}+(-1)^{r-1}(-1)^{p-1}\det J_{p-2} \cdot \det \overline{J_{q-p-1}} \cdot \det J_{r-q-1}^{\prime\prime} \cdot \det \tilde{J}_{r-1}^{\prime\prime} \\
&\hspace{2cm}-(-1)^{r-1}(-1)^q2\det J_{p-2} \cdot \det J_{q-p-1}\cdot \det \overline{J_{r-q-1}} \cdot \det \tilde{J}_{r-1}^{\prime\prime} \\
&+2(-1)^{q-1}(-1)^{p-1}\det J_{p-2} \cdot \det \overline{J_{q-p-1}'} \cdot \det J_{r-q-1}\cdot \det \tilde{J}_{r-1}\\ 
&\hspace{1cm}+2(-1)^q(-1)^{p-1}\det J_{p-2} \cdot \det \overline{J_{q-p-1}} \cdot \det J_{r-q-1}^{\prime\prime}  \cdot \det \tilde{J}_{r-1}\\
&\hspace{2cm}+4\det J_{p-2} \cdot \det J_{q-p-1}\cdot \det J_{r-q-1}\cdot \det \tilde{J}_{r-1}\\
 &=(-1)^{r-1}r(p-1).
\end{align*}
Finally
{\small
\begin{equation*}
\cofac_{23}=\left| \;
\begin{array}{cc|cc|cc|cc|cc}
  &&&&&&&&& \\ 
  J_{p-2}&&&&&&&&&\\
  &&&&&&&&1&\\
  \hline
  &&&&&&&&1&\\
  &&J_{q-p-1}&&&&&&&\\
  &&&&&&&&&1\\
  \hline
  &&&&&&&&&1\\
  &&&&J_{r-q-1}&&&&\\
  &&&&&&&&&\\
   \hline
  &&&&&&&&&\\
  &&&&&&&\tilde{J}_{r-1}&&\\
  &&&&&&&&&\\
  \hline
  &1&1&&&&&&2&\\
  &&&&&1&-1&&&0\\
 \end{array}
  \; \right|
\end{equation*}
} 
which can be calculated with respect to the last two rows so that
\begin{align*}
\cofac_{23}
&=(-1)^{r-1}(-1)^{q-1}\det J_{p-2}^{\prime}\cdot\det J_{q-p-1}\cdot\det\overline{ J_{r-q-1}^{\prime}}\cdot\det \tilde{J}_{r-1}\\ 
&\hspace{1cm}+(-1)^{r-1}(-1)^{q-1}\det J_{p-2}\cdot\det J_{q-p-1}^{\prime\prime}\cdot\det \overline{J_{r-q-1}^{\prime}}\cdot\det \tilde{J}_{r-1} \\
&\hspace{2cm}-(-1)^{r-1}(-1)^{q-1}2\det J_{p-2}\cdot\det J_{q-p-1}\cdot\det \overline{J_{r-q-1}^{\prime}}\cdot\det \tilde{J}_{r-1} \\
&+(-1)^{r-1}(-1)^{q-1}\det J_{p-2}^{\prime} \cdot \det J_{q-p-1} \cdot \det \overline{J_{r-q-1}} \cdot \det \tilde{J}_{r-1}^{\prime\prime}\\
&\hspace{1cm}+(-1)^{r-1}(-1)^{q-1}\det J_{p-2} \cdot \det J_{q-p-1}^{\prime\prime}\cdot \det \overline{J_{r-q-1}} \cdot \det \tilde{J}_{r-1}^{\prime\prime}\\
&\hspace{2cm}-(-1)^{r-1}(-1)^{q-1}2\det J_{p-2} \cdot \det J_{q-p-1}\cdot \det \overline{J_{r-q-1}} \cdot \det \tilde{J}_{r-1}^{\prime\prime}\\
 &=-(-1)^{r-1}r(q-1).
\end{align*}
As a result, we conclude that

\begin{align*}
 c^2(W)&=(-2u,-2v,-2y) \left(\begin{array}{ccc}
p(p-1)&q(p-1)&r(p-1)\\
q(p-1)&q(q-1)&r(q-1)\\
r(p-1)&r(q-1)&r(r-1)\\
\end{array}\right)(-2u,-2v,-2y)^{T}\\
&=4u^2p(p-1)+4v^2q(q-1)+4y^2r(r-1)\\
&\quad+8uvq(p-1)+8uyr(p-1)+8vyr(q-1).
\end{align*}

\end{document}